\documentclass{lmcs}
\pdfoutput=1
\usepackage[utf8]{inputenc}

\usepackage{lastpage}
\lmcsdoi{21}{2}{22}
\lmcsheading{}{\pageref{LastPage}}{}{}%
{May~09,~2024}{Jun.~10,~2025}{}

\keywords{indivisibility, Weihrauch reducibility, Ramsey theory, reverse
mathematics, computable structure theory, computable combinatorics}

\usepackage{hyperref}

\usepackage{rotating}
\usepackage{amssymb}
\usepackage{gensymb}
\usepackage{amsfonts}
\usepackage{eucal}
\usepackage{theoremref}
\usepackage{mathtools}
\usepackage{mathrsfs}
\usepackage{graphicx}
\usepackage{changepage}
\usepackage{adjustbox}
\usepackage{microtype}

\usepackage{tikz,tikz-cd,pgf,pgfplots}
\pgfplotsset{compat=1.18}

\newcommand{\N}{\mathbb{N}}
\newcommand{\Q}{\mathbb{Q}}
\newcommand{\calA}{\mathcal{A}}
\newcommand{\calS}{\mathcal{S}}

\DeclareMathOperator{\dom}{\mathrm{dom}}
\DeclareMathOperator{\ran}{\mathrm{ran}}
\DeclareMathOperator{\Det}{\mathrm{Det}}

\newcommand{\restr}{\!\upharpoonright\!}
\newcommand{\abs}[1]{\left\lvert{#1}\right\rvert}

\newcommand{\cvg}{\mkern-5mu\downarrow}

\newcommand{\bdef}[1]{\emph{#1}}

\newcommand{\W}{\mathrm{W}}
\newcommand{\sW}{\mathrm{sW}}
\newcommand{\rmc}{\mathrm{c}}
\newcommand{\rmsc}{\mathrm{sc}}

\newcommand{\RT}{\mathsf{RT}}
\newcommand{\cRT}{\mathsf{cRT}}
\newcommand{\SRT}{\mathsf{SRT}}
\newcommand{\COH}{\mathsf{COH}}
\newcommand{\CN}{\mathsf{C}_\N}
\newcommand{\TCN}{\mathsf{TC}_\N}
\newcommand{\Eqinf}{\mathscr{E}}
\newcommand{\FE}{\mathrm{FE}}

\DeclarePairedDelimiterX\set[1]\lbrace\rbrace{\,\def\st{\,:\,}#1\,}
\DeclarePairedDelimiterX\ang[1]\langle\rangle{\mkern-.5mu\def\st{\,:\,}#1\mkern-.5mu}

\DeclareMathOperator{\Indtext}{\mathsf{Ind}}
\newcommand{\Ind}[2][]{\Indtext{#2}_{{#1}}}

\begin{document}

\title{Indivisibility and uniform computational strength}
\thanks{This work extends results from the first chapter of the author's
    Ph.D.\ dissertation at Penn State University \cite{mythesis}. The author is
    grateful to his thesis advisors Linda Westrick and Jan Reimann for their
    invaluable help and support, to Arno Pauly for several insightful
    conversations surrounding $\Ind{\Q}$, and to Isabella Scott and Reed Solomon
    for pointing out errors in earlier versions of this paper. He would also
    like to thank the anonymous referee for comments leading to a much improved
    presentation in several places.
    This research was supported in part by NSF grant DMS-1854107.}

\author[K.~Gill]{Kenneth Gill\lmcsorcid{0009-0002-3452-3937}}

\address{Department of Mathematics and Computer Science\\ La Salle University\\
Philadelphia, PA, USA}
\email{gillmathpsu@posteo.net}


\begin{abstract}
  A countable structure is indivisible if for every coloring with finite range there is a monochromatic isomorphic subcopy of the structure. Each indivisible structure naturally corresponds to an indivisibility problem which outputs such a subcopy given a presentation and coloring. We investigate the Weihrauch complexity of the indivisibility problems for two structures: the rational numbers $\mathbb{Q}$ as a linear order, and the equivalence relation $\mathscr{E}$ with countably many equivalence classes each having countably many members. We separate the Weihrauch degrees of both corresponding indivisibility problems from several benchmarks, showing in particular that the indivisibility problem for $\mathbb{Q}$ cannot solve the problem of finding a monochromatic rational interval given a coloring for which there is one; and that the Weihrauch degree of the indivisibility problem for $\mathscr{E}$ is strictly between those of $\mathsf{RT}^2$ and $\mathsf{SRT}^2$, two widely studied variants of Ramsey's theorem for pairs whose reverse-mathematical separation was open until recently.
\end{abstract}

\maketitle


\section{Introduction}

Reverse mathematics is a branch of logic which seeks to classify theorems based
on their logical strength. This is often done by working over a relatively
weak base system for mathematics, such as a fragment of second-order arithmetic,
and then looking at the truth of various implications in models of that fragment.
It turns out that one can frequently obtain more fine-grained distinctions
between theorems by comparing their intrinsic computational or combinatorial
power, rather than being limited to considerations of provability only. For
theorems expressible as $\Pi^1_2$ statements of second-order arithmetic, this
can be achieved by identifying a theorem with an instance-solution pair, where
instances are all sets $X$ to which the hypothesis of the theorem apply, and
solutions to $X$ are all sets $Y$ witnessing the truth of the conclusion of the
theorem for $X$.
One then introduces various reducibilities to compare, among different
instance-solution pairs, the difficulty of finding a solution from a given
instance.

Instance-solution pairs can be formalized as problems:

\begin{defi}
A \bdef{problem} is a partial multivalued function on Baire space, $P\colon
\subseteq \N^\N \rightrightarrows \N^\N$. Any $x\in \dom P$ is called an
\bdef{instance} of $P$, and any $y\in P(x)$ is a \bdef{solution} to $x$ with
respect to $P$. 
\end{defi}

In the present work, we focus on problems arising from indivisible structures,
and use mainly Weihrauch reducibility to gauge their uniform computational content.
We assume familiarity with computability theory and the basic terminology of
computable structure theory, as in the introductory sections of \cite{Mon1}.
All our structures are countable and we will take them to have domain $\N$
when convenient.

\begin{defi}
A structure $\calS$ is \textit{indivisible} if for any presentation
$\calA$ of $\calS$ and any coloring of $A = \dom \calA$ with finitely many
colors, there is a monochromatic subset $B\subseteq A$ such that $\mathcal{B}
\simeq \calA$, where $\mathcal{B}$ is the substructure induced by $B$.
$\calS$ is \bdef{computably indivisible} if in addition such a set $B$ can
always be computed from $\calA$ and $c$.
\end{defi}

Any single-element structure is indivisible, but otherwise an indivisible
structure must be infinite, and we will take all structures to be infinite
unless otherwise specified.
Indivisibility belongs properly to Ramsey theory, the study of combinatorial
structures which cannot be ``too disordered'' in that any large subset must
contain a highly organized substructure. This field goes back to 1930, when
Ramsey's theorem first appeared, and principles arising in Ramsey theory have
long been objects of interest in reverse mathematics.

\begin{defi} If $\calS$ is an indivisible structure, then $\Ind{\calS}$
  is the \bdef{indivisibility problem} associated to $\calS$, for which
\begin{itemize}
  \item Instances are triples $\ang{\calA,c,k}$ where $\calA$ is a presentation of
    $\calS$ (identified with its atomic diagram) and $c$ is a $k$-coloring of
    $A$, and
  \item Solutions to $\ang{\calA,c,k}$ are (characteristic functions of) subsets
    $B$ of $A$ which are monochromatic for $c$ and such that $\mathcal{B} \simeq
    \calA$.
\end{itemize}

Here $\calA$ is allowed to be an arbitrary presentation of $\calS$, not
necessarily one having domain $A = \N$. We still treat $c$ as a coloring of all
of $\N$ but disregard any color assigned to a point outside $A$. By convention,
$i \in A$ iff the formula $x_i = x_i$ is true in the atomic diagram of $\calA$,
and if it is false then all atomic formulas involving $x_i$ should also be set
to false.
We will also use the following variations on $\Ind{\calS}$:
\begin{itemize}
  \item For each fixed $k$, let $\Ind[k]{\calS}$ have domain $\set{ \ang{\calA,c}
    \st \ang{\calA,c,k} \in \dom (\Ind{\calS})}$ and let $\Ind[k]{\calS}(
    \ang{\calA,c}) = \Ind{\calS}(\ang{\calA,c,k})$. 
    That is, $\Ind[k]{\calS}$ is the restriction of $\Ind{\calS}$ to $k$-colorings,
    so that $\Ind{\calS}$ is exactly the problem $\bigsqcup_{k\in \N}
    \Ind[k]{\calS}$, where $\sqcup$ is defined in the next section.
  \item Let $\Ind[\N]{\calS}$ be the same problem as $\Ind{\calS}$ except
    without any given information on the number of colors used.
    In other words, $\dom (\Ind[\N]{\calS}) =
    \bigcup_{k\in \N} \dom (\Ind[k]{\calS})$ and $\Ind[\N]{\calS}(x) =
    \Ind[k]{\calS}(x)$ if $x\in \dom (\Ind[k]{\calS})$. 
\end{itemize}
\end{defi}
One can also define the \bdef{strong indivisibility problem} where a solution
additionally includes an explicit isomorphism between $\mathcal{B}$ and $\calA$.
We do not study this notion here; both it and $\Ind{\calS}$, or more
specifically what we call $\Ind[2]{\calS}$, are discussed in unpublished work of
Ackerman, Freer, Reimann, and Westrick \cite{Ind}.

This paper focuses on two indivisible structures in particular. First, the set
of rational numbers $\Q$ is viewed as, up to order isomorphism, the unique
countable dense linear order with no greatest or least element. $\Ind{\Q}$ is
studied in Section~\ref{sec:rats}. The other structure is the so-called
countable equivalence relation:

\begin{defi} The \bdef{countable equivalence relation} $\Eqinf$ is, up to
  isomorphism, the structure in the language $\{R\}$ of a single binary relation
  such that $R$ is an equivalence relation, and such that $E$ is divided into
  countably many equivalence classes each having countably many members.
\end{defi}

To see $\Eqinf$ is indivisible, suppose its elements are colored red and blue
(it suffices to show for two colors). Then either there are infinitely many
equivalence classes each containing infinitely many red points, or there are
cofinitely many equivalence classes each containing cofinitely many blue points.
Either way one gets a monochromatic subcopy of $\Eqinf$.
$\Ind{\Eqinf}$ is studied in Section~\ref{sec:eq}. Our main results separate
$\Ind{\Q}$ and $\Ind{\Eqinf}$ from various benchmark problems.
In the case of $\Q$, we have
\begin{itemize}
  \item $\Ind{\Q}$ is Weihrauch incomparable with $\CN$ (\thref{QperpCN}), which
    can be viewed as the problem which finds a monochromatic rational interval
    given a coloring for which there is one.
  \item $\Ind[k]{\Q}$ cannot be solved by any problem having a c.e.\ approximation
    (\thref{pcet}), in a precise sense given by \thref{pcetdef} which
    encompasses $\RT^1_\N$ and any problem with finite computable codomain,
    among other things.
\end{itemize}

The hierarchy of Weihrauch reductions (and nonreductions) relating
versions of $\TCN^k$, $\Ind[k]{\Q}$, and $\RT^1_k$ is shown in Figure~\ref{Qreds}
below. 
Next, $\Ind{\Eqinf}$ turns out to be closely related to but distinct from Ramsey's
theorem for pairs:
\begin{itemize}
  \item The Weihrauch degree of $\Ind[k]{\Eqinf}$ is strictly between those of
    $\SRT^2_k$ and $\RT^2_k$ for all $k\geq 2$;
  \item In particular, $\Ind[2]{\Eqinf}$ is not Weihrauch reducible to $\SRT^2_\N$
    and $\RT^2_2$ is not computably reducible to $\Ind[\N]{\Eqinf}$ (\thref{eqthm}).
\end{itemize}
The Weihrauch reductions between these principles are displayed in
Figure~\ref{eqreds}. 

Before getting to our main results in Sections~\ref{sec:rats} and \ref{sec:eq},
we define all the reducibilities under consideration in Section~\ref{sec:defs}
as well as the main problems used as benchmarks. In Section~\ref{sec:rmk} we
comment on the role of uniform computable categoricity in the present
investigation and collect several properties of indivisibility problems in
general. Finally, in Section~\ref{sec:moredirs} we briefly mention some
directions for future research on $\Ind{\Q}$ and $\Ind{\Eqinf}$.


\section{Preliminaries}\label{sec:defs}

The standard reference on reverse mathematics is Simpson \cite{Simpson}. 
The books by Hirschfeldt \cite{STT} and by Dzhafarov and Mummert \cite{DM} are
more computability-oriented treatments of the subject, the latter discussing in
detail the reducibilities we employ. 
Our notation and conventions with respect to computability theory are standard
and can be found in any modern textbook, e.g., Soare \cite{Soare}.
For any set $X$ and $n\in \N$, $X\restr n$ is the truncation of $X$ up through
its first $n$ bits.
If $\sigma\in \N^{<\N}$, then $[\sigma]$ denotes the set of finite or infinite
strings extending $\sigma$.
We will always identify $k\in \N$ with the set $\{0,\dotsc, k-1\}$. Write $W_e$
for the $e$th c.e.\ set with respect to some computable enumeration of Turing
machines, and $W_e^X$ for the $e$th $X$-c.e.\ set.

\begin{defi}
Let $P$ and $Q$ be problems. Then 
\begin{itemize}
  \item $P$ is \bdef{computably reducible} to $Q$, written $P\leq_\rmc Q$, if
    each $p\in \dom P$ computes some $x\in \dom Q$ such that for every $q\in
    Q(x)$, $p\oplus q$ computes an element of $P(p)$. If we require $q$ by
    itself to compute an element of $P(p)$, then $P$ is \bdef{strongly
    computably reducible} to $Q$, written $P \leq_\rmsc Q$.
  \item $P$ is \bdef{Weihrauch reducible} to $Q$, written $P\leq_\W Q$, if there
    are Turing functionals $\Delta$ and $\Psi$ such that if $p\in \dom P$, then
    $\Delta^p\in \dom Q$; and for any $q\in Q(\Delta^p)$, we have $\Psi^{p\oplus
    q}\in P(p)$. If we require $\Psi$ to only have oracle access to $q$, then
    $P$ is \bdef{strongly Weihrauch reducible} to $Q$, written $P \leq_\sW Q$.
\end{itemize}
\end{defi}

$\Delta$ and $\Psi$ are sometimes called the forward and return functionals,
respectively.
We have $P \leq_\sW Q \implies P \leq_\W Q \implies P \leq_\rmc Q$ and $P \leq_\sW
Q \implies P \leq_\rmsc Q \implies P \leq_\rmc Q$, but none of these implications
reverse and there is no logical relationship between $\leq_\W$ and $\leq_\rmsc$.
All of these reducibilities also imply that $Q$ implies $P$ in $\omega$-models
of the weak logical system $\mathsf{RCA}_0$, i.e., models with the standard
natural numbers. 
Write $P\vert_\square Q$ if $P$ and $Q$ are incomparable under $\leq_\square$,
and $P \equiv_\square Q$ if $P\leq_\square Q$ and $Q\leq_\square P$, for
$\square \in \{ \rmc,\rmsc,\W,\sW \}$. In the latter case $P$ and $Q$ are said
to be equivalent or to have the same degree with respect to $\leq_\square$.
See \cite{HJ16,DM} for details on these and other reducibilities, and \cite{BGP} for a
good recent survey of Weihrauch reducibility with a large bibliography and
historical remarks. The reducibilities above can be defined between problems on
any represented spaces, but we do not need this level of generality. More
information can again be found in \cite{BGP}, and we largely follow the latter
treatment here.

There are a few algebraic operations on problems which we will have occasion to
use. Let $P$ and $Q$ be problems. Then
\begin{itemize}
  \item $P\times Q$ is the \bdef{parallel product} of $P$ and $Q$, defined by
    $\dom(P\times Q) = \dom P \times \dom Q$, and $(P\times Q)(p,q) = P(p)
    \times Q(q)$.
  \item $P^n$ is the $n$-fold parallel product of $P$ with itself.
  \item $P^\ast$ is the \bdef{finite parallelization} of $P$, which solves $P^n$
    for any $n$. Instances include the data \nolinebreak $n$.
  \item $\hat{P}$ is the \bdef{parallelization} of $P$. Instances are sequences
    $\ang{p_0,p_1,p_2,\dotsc} \in \dom(P)^\N$, and the set of solutions to this
    instance is the Cartesian product $\displaystyle \prod_{i\in \N} P(p_i)$.
  \item $P \sqcup Q$ is the \bdef{coproduct} of $P$ and $Q$, with $\dom(P \sqcup
    Q) = \dom P \sqcup \dom Q$, and $(P\sqcup Q)(0,p) = \{0\} \times P(p)$ and
    $(P\sqcup Q)(1,q) = \{1\} \times Q(q)$. So $P\sqcup Q$ is the problem which
    is capable of solving any instance of $P$ or of $Q$, but only one at a time.
  \item $P'$ is the \bdef{jump} of $P$. An instance of $P'$ is a sequence $\ang{p_0,
    p_1, p_2,\dotsc}$ of reals converging to some $p\in \dom P$ in the Baire
    space topology (equivalently, converging entrywise). The solutions to this
    instance are just the elements of $P(p)$. So $P'$ answers the same question
    about $p\in \dom P$ that $P$ does, but only has access to a limit
    representation of $p$. We also set $P^{(n+1)} = (P^{(n)})'$ with $P^{(0)} =
    P$.
  \item $P\ast Q$ is the \bdef{compositional product} of $P$ and $Q$. This can
    be characterized intuitively as the strongest problem under $\leq_\W$
    obtainable as the composition of $f$ and $g$, ranging over all $f\leq_\W P$
    and $g\leq_\W Q$.
    (A more precise definition of a representative of the degree of $P\ast Q$ is
    given in Section~\ref{sec:rmk} below.)
\end{itemize}

The following ``benchmark'' problems will be used as a basis for comparison with
the problems we study:
\begin{itemize}
  \item $\mathsf{LPO}$, the \bdef{limited principle of omniscience}, has
    $\dom \mathsf{LPO} = \N^\N$, with $\mathsf{LPO}(0^\N) = 0$ and
    $\mathsf{LPO}(p) =1$ otherwise.
    $\mathsf{LPO}^{(n)}$ can be thought of as answering a single
    $\Sigma^0_{n+1}$ question.
  \item $\lim$ maps a convergent sequence of reals $\ang{p_0,p_1,p_2,\dotsc}$ to
    its limit. 
  \item $\CN$, closed choice on $\N$, outputs an element of a
    nonempty set $A\subseteq\N$ given an enumeration of its complement.
  \item $\TCN$, the totalization of $\CN$, extends $\CN$ by allowing $A =
    \emptyset$ and outputting any number in this case.
  \item $\RT^n_k$, Ramsey's theorem for $n$-tuples and $k$ colors, has instances
    $c\colon [\N]^n\to k$, and solutions to $c$ are (characteristic functions of)
    infinite $c$-homogeneous sets. Here $[X]^n$ is the set of $n$-element
    subsets of $X$, and a set $X\subseteq \N$ is \bdef{homogeneous} for $c$ if
    $X$ is infinite and $c$ is monochromatic on $[X]^n$.
  \item $\SRT^2_k$, stable Ramsey's theorem for pairs, is the restriction of
    $\RT^2_k$ to \bdef{stable colorings}, i.e., colorings $c$ such that $\lim_m
    c\{n,m\}$ exists for all $n$.
  \item $\RT^n_\N$ has instances $c\in \bigcup_{k\in\N} \dom
    \RT^n_k$. Solutions to $c$ are again infinite $c$-homogeneous sets. Notice
    that in this formulation, the number of colors used is not included as part
    of an instance. $\SRT^2_\N$ is defined similarly.
  \item For $k\in \N \cup \{ \N \}$, $\cRT^n_k$ is the ``color version'' of
    $\RT^n_k$, which only outputs the colors of $\RT^n_k$-solutions. We have
    $\RT^1_k \equiv_\W \cRT^1_k$ for all $k$, since the color can be used to
    compute the set of points of that color and vice versa.
  \item $\RT^n_+$ is $\bigsqcup_{k\in \N} \RT^n_k$, and $\cRT^n_+$ and
    $\SRT^2_+$ are defined similarly.
\end{itemize}

Clearly $\RT^n_2 \leq_\sW \RT^n_3 \leq_\sW \dotsc \leq_\sW \RT^n_+ \leq_\sW
\RT^n_\N$ for all $n$, and similarly for $\SRT^2_k$. Also $\Ind[2]{\calS}
\leq_\sW \Ind[3]{\calS} \leq_\sW \dotsc \leq_\sW \Ind{\calS} \leq_\sW
\Ind[\N]{\calS}$ for any indivisible $\calS$.


\section{Uniform computable categoricity}\label{sec:rmk}

Recall that a computable structure $\calA$ is uniformly computably categorical
(u.c.c.)
if there is a Turing functional which computes an isomorphism from $\mathcal{B}$
to $\calA$ given the atomic diagram of any presentation $\mathcal{B}$ of
$\calA$. If $\calS$ is indivisible and u.c.c., then
from the point of view of $\leq_\W$ we can as a convention regard
the instances of $\Ind{\calS}$ as only including the data $k$ and $c\colon \N \to
k$, since in a reduction, the functionals $\Delta$ and $\Psi$ can just build in
the translations between any given presentation $\calA$ and some fixed
computable presentation of $\calS$. Instances of $\Ind[k]{\calS}$ and
$\Ind[\N]{\calS}$ are viewed simply as colorings $c$ of $\N$, which is justified
as long as we choose a computable reference presentation $\calS$. 
One can state this more formally as

\begin{prop}\thlabel{ignorepres} If $\calS$ is indivisible and u.c.c.,
  then $\Ind[k]{\calS} \equiv_\W P_k$ where $P_k$ is the
  restriction of $\Ind[k]{\calS}$ to instances of the form $\ang{\calS, c}$.
  Analogous statements hold for $\Ind{\calS}$ and $\Ind[\N]{\calS}$.\qed
\end{prop}

These conventions are not a priori applicable when considering $\leq_\sW$ (or
$\leq_\rmsc$),
because the return functional $\Psi$ could need oracle access to
the presentation $\calA$ in order to translate back to a solution of
$\ang{\calA, c, k}$ if $\Delta$ modified $\calA$. A version of
\thref{ignorepres} with $\equiv_\W$ replaced by $\equiv_\rmc$ still holds if
$\calS$ is merely relatively computably categorical, i.e., if the functional
computing the isomorphism from $\mathcal{B}$ to $\calA$ (as above) is allowed to
depend on $\mathcal{B}$.\footnote{Technically, what we call ``uniform computable
  categoricity'' should really be called uniform relative computable
  categoricity. The standard definition of uniform computable categoricity of a
  computable structure $\calA$ is that there is a Turing functional
  which, given any \emph{computable} copy $\mathcal{B}$ of $\calA$ as oracle,
  computes an isomorphism from $\mathcal{B}$ to $\calA$. However, this is
  equivalent to our definition by a theorem of Ventsov
\cite[Theorem~III.18]{Mon1}. Relative and plain computable categoricity differ;
see Chapter~VIII of \cite{Mon1} for details.}
\thref{IndSstar} below also continues to hold for relatively computably
categorical $\calS$, so long as $\equiv_\W$ is changed to $\equiv_\rmc$
everywhere.

Both of the structures we focus on in this paper are uniformly computably
categorical and so we will follow the above convention without further comment.
That $\Q$ is u.c.c.\ follows from the classical back-and-forth argument, which
is effective.
To see that $\Eqinf$ is u.c.c., for any presentation
$\calA$ of $\Eqinf$, decompose any $n\in A$ as a pair $\ang{x,y}_\calA$ so that
$n$ is the $y$th element of the $x$th distinct equivalence class, in order of
discovery within the atomic diagram of $\calA$. Then if also
$m=\ang{z,w}_\calA$, we have that $n$ and $m$ are equivalent iff $x=z$, and
$\ang{\cdot,\cdot}_\calA$ is uniformly $\calA$-computable. (This definition does
not uniquely specify $\ang{\cdot,\cdot}_\calA$, but one can make some canonical
choice.)
If we take $\Eqinf$ to be computable and $\calA$ is a given copy of $\Eqinf$,
without loss of generality with $E=A=\N$, then the map $\ang{x,y}_\calA \mapsto
\ang{x,y}_\Eqinf$ is a uniformly $\calA$-computable isomorphism between $\calA$
and $\Eqinf$.

The problems corresponding to u.c.c.\ indivisible structures enjoy several nice
properties. For example, it is not hard to see that the indivisibility and
strong indivisibility problems of such a structure must be Weihrauch equivalent.
There is also the following easy observation:

\begin{prop}\thlabel{IndSstar} If $\calS$ is indivisible and u.c.c., then
  $\Ind{\calS} \equiv_\W \Ind{\calS}^\ast$. Similarly for $\Ind[\N]{\calS}$.
\end{prop}
\begin{proof}
  It suffices to show $\Ind{\calS}^\ast \leq_\W \Ind{\calS}$. 
  Since $\calS$ is u.c.c., as noted above, we can
  assume each instance of $\Ind{\calS}$ consists only of a coloring together
  with a number of colors. 
  Then suppose we are given instances $\ang{c_0,k_0}, \dotsc, \ang{c_n,k_n}$ of
  $\Ind{\calS}$ in parallel, with $c_i \colon \N \to k_i$ for each $i$. 
  Let $p_i$ be the $i$th prime number and define $d\colon \N \to \prod_i
  p_i^{k_i+1}$ by $d(x) = \prod_i p_i^{c_i(x)+1}$. If $H$ is a solution to $d$
  with color $j$, then $j = \prod_i p_i^{j_i+1}$ for some numbers $j_0 < k_0,
  \dotsc, j_n < k_n$. Hence $c_i(H) = j_i$ for each $i=0, \dotsc, n$, and $H$ is
  simultaneously a solution to $\ang{c_i,k_i}$ for all $i$.
  It is clear that the argument works unchanged for $\Ind[\N]{\calS}$.
\end{proof}

Next we have the following result which extends Proposition~62 of
\cite{PPS23} from $\N$ to any indivisible, u.c.c.\ $\calS$. Our argument is
really only a slight adaptation of the original, but we give it in detail for
completeness. It will be useful for the proof to use the representative $f\star
g$ of the Weihrauch degree of $f\ast g$, defined by
\[
  (f\star g)(x,y) = \ang{\mathrm{id} \times f} \circ \Phi_x \circ g(y),
\]
where $\mathrm{id}$ is the identity map on $\N^\N$ and $\Phi$ is a universal
functional. (See \cite[Definition~11.5.3]{BGP}.)

\begin{prop}\thlabel{IndSN} For any indivisible $\calS$, we have $\Ind[\N]{\calS}
  \leq_\sW \Ind{\calS} \ast \CN$. If $\calS$ is u.c.c.,
  then also $\Ind{\calS} \ast \CN \leq_\W \Ind[\N]{\calS}$.
\end{prop}
\begin{proof}[Proof, after \cite{PPS23}] 
  For the first statement, let $\calS$ be an arbitrary indivisible structure and
  let $\ang{\calA,c}$ be a given instance of $\Ind[\N]{\calS}$. Build a
  $\CN$-instance by, whenever a new color is seen in $c$, enumerating all
  numbers less than that color. One can use a $\CN$-solution $k$ together with
  $\ang{\calA,c}$ (which can be encoded as part of the program to be run by
  $\Phi$) to compute the instance $\ang{\calA,c,k}$ of $\Ind{\calS}$, and a
  solution to this instance is also a solution to $\ang{\calA,c}$.

  To prove the second statement, let $\mathsf{Bound}$ be the problem that
  outputs an upper bound on an enumeration of a finite set. This is Weihrauch
  equivalent to $\CN$ and it will be convenient to show that $\Ind{\calS} \star
  \mathsf{Bound} \leq_\W \Ind[\N]{\calS}$. Let $(x,Y)$ be an instance of $\Ind{\calS}
  \star \mathsf{Bound}$, so $Y$ is a finite set represented by an (infinite)
  enumeration. Since $\calS$ is u.c.c., we can treat $\Phi_x$ in its second
  component as computing only a coloring with an upper bound on its range. 
  For each $i\in \N$, let $p_i$ be the $i$th prime number as before. Build $d
  \in \dom (\Ind[\N]{\calS})$ as follows: find a number $i_0$ such that
  $\Phi_x(i_0)$ outputs, in addition to a partial coloring $c_0$ of $\N = \dom
  \calS$, an upper bound $k_0$ on the range of $c_0$. This $i_0$ is an initial
  guess for an element of $\mathsf{Bound}(Y)$. Set $d(n)=p_{i_0}^{c_0(n)}$ for
  all $n$ for which we see $c_0(n)$ defined by $\Phi_x(i_0)$. In general, if
  $i_s$ has been found, let $d(n) = p_{i_s}^{c_s(n)}$ whenever $c_s(n)$ is
  defined and $n$ had not been colored at a previous stage. If numbers are
  enumerated into $Y$ so that $\max Y \geq i_s$, find an $i_{s+1} > \max Y$ such
  that $\Phi_x(i_{s+1})$ outputs a number $k_{s+1}$; such an $i_{s+1}$ must
  exist, so we can hold off on extending $d$ until $k_{s+1}$ appears. When it
  does, continue with $d(n) = p_{i_{s+1}}^{c_{s+1}(n)}$ for all $n$ which had
  not been previously colored and for which $c_{s+1}(n)$ is defined.

  Eventually this process stabilizes at some $i_\ell$, since $Y$ is finite, and
  once it does $\Phi_x(i_\ell)$ must produce a total coloring. Then any solution
  $H$ of $d$ has color $p_{i_\ell}^a$ for some $a < k_\ell$. By rerunning the
  procedure in the last paragraph, the return functional can recover $d$, hence
  find the color of $H$, and from that learn $i_\ell$ and output the first
  component of $\Phi_x(i_\ell)$ to satisfy $\mathrm{id}$. Finally, $H$ is in
  fact a solution to the second component of $\Phi_x(i_\ell)$, because it can
  only include points colored after $i_\ell$ stabilizes and so the fact that $d$
  and the coloring computed by $\Phi_x(i_\ell)$ differ on finitely many other
  points is of no consequence. \qedhere
\end{proof}

If $\calS$ is not u.c.c., then $d$ could still be built to eventually agree
(up to taking prime powers) with the correct coloring $c_\ell$, so that it has
the same solutions. 
But one would also need to somehow copy the atomic diagram computed by
$\Phi_x(i_s)$ to produce an instance of $\Ind[\N]{\calS}$, and there is no
reason something eventually agreeing with a correct presentation should be a
presentation itself at all, never mind one for which $d$ has the same
solutions.
On the other hand, exactly the same proof serves to show that $\RT^n_\N
\equiv_\W \RT^n_+ \ast \CN$ for all $n$ and that $\SRT^2_\N \equiv_\W \SRT^2_+
\ast \CN$.

The nonuniform situation is considerably simpler. 
$P\ast \CN \equiv_\rmc P$ for any problem $P$ since the solutions to $\CN$ are
just single numbers, so the problems $\Ind{\calS}$, $\Ind[\N]{\calS}$, and
$\Ind{\calS} \ast \CN$ collapse under $\leq_\rmc$, and indeed
\begin{prop}\thlabel{IndSNsc} For any indivisible $\calS$, we have
  $\Ind{\calS} \equiv_\rmsc \Ind[\N]{\calS}$.
\end{prop}
\begin{proof} $\Ind{\calS} \leq_\rmsc \Ind[\N]{\calS}$ is always true, and
  $\Ind[\N]{\calS} \leq_\rmsc \Ind{\calS}$ follows because one can nonuniformly
  append the numbers of colors used by an $\Ind[\N]{\calS}$-instance in order to
  get an $\Ind{\calS}$-instance with the same solutions.
\end{proof}

We now discuss one more general property which will be useful later. 
A problem $f$
is called a \bdef{fractal} if there is an $F\colon \subseteq \N^\N
\rightrightarrows \N^\N$ with $f\equiv_\W F$ such that $F \restr_A \equiv_\W F$
for every clopen $A$ with $A\cap \dom F \neq \emptyset$, and $f$ is a
\bdef{total} or \bdef{closed fractal} if $F$ can additionally be taken with
domain $\N^\N$ \cite[Definition~11.4.10]{BGP}. Intuitively, being a fractal
means that no matter how far one zooms into the domain of $f$, its full
power is always present.

We do not expect $\Ind[k]{\calS}$ to be a fractal even for an arbitrary u.c.c.\
structure, but it turns out to be one for both $\calS = \Q$ and $\Eqinf$ because
both structures have the property that deleting any finite set produces an
isomorphic substructure. 
This observation can be made more general by invoking the notion
of finite tolerance, originally introduced in \cite{Dorais} (see
\cite[Definition~11.4.5]{BGP}). A problem $f$ is \bdef{finitely tolerant} if
there is a computable functional $\Phi$ such that for any $x\in \dom f$, for any
$y\in f(x)$, we have that $\Phi(y,n)$ computes a solution to $f(z)$ for any
instance $z$ such that $z(i)=x(i)$ for all $i\geq n$.
Thus $\Ind[k]{\Q}$ and $\Ind[k]{\Eqinf}$ are finitely tolerant: in each case,
given a solution $H$ to a coloring $c$ along with $n$, $\Phi$ can just delete
the first $n$ elements from $H$ to produce a solution to any coloring $d$
agreeing with $c$ after the first $n$ bits. The same trick applies to any
structure which is isomorphic to every cofinite substructure.
(Not every indivisible structure has the latter property; consider $(\Q+1,<)$.)

\begin{prop}\thlabel{SNfractal} ~
  \begin{enumerate}[label=(\roman*)]
    \item If $\calS$ is u.c.c.\ and $\Ind[k]{\calS}$ is finitely tolerant, then
      $\Ind[k]{\calS}$ is a total fractal.
    \item $\Ind[\N]{\calS}$ is a fractal for any u.c.c.\ $\calS$.
  \end{enumerate}
\end{prop}
\begin{proof} ~
  \begin{enumerate}[label=(\roman*)]
    \item If $x$ is any real, define $\lfloor x\rfloor_k$ by setting $\lfloor
      x\rfloor_k(i)=\min\{k-1,x(i)\}$ for each $i$. Then let $f \colon \N^\N
      \rightrightarrows \N^\N$ be given by $f(x) = \Ind[k]{\calS}(\lfloor
      x\rfloor_k)$. Viewing the instances of $\Ind[k]{\calS}$ as colorings of a
      fixed computable presentation, it is clear that $f \equiv_\W
      \Ind[k]{\calS}$: on the one hand, $\Ind[k]{\calS}$ is a subproblem of $f$,
      and on the other hand, we can reduce $f$ to $\Ind[k]{\calS}$ by computing
      $\lfloor x\rfloor_k$ from a given $x$ and feeding it into
      $\Ind[k]{\calS}$.  Therefore the latter is Weihrauch equivalent to a total
      problem. Notice $f(x) = f(\lfloor x\rfloor_k)$ for all $x$. (So far we
      have not used the assumption of finite tolerance.) 

      We now show that this $f$ is a fractal by reducing it to its restriction
      to a given cylinder. Let $\Phi$ witness finite tolerance of
      $\Ind[k]{\calS}$. Given a finite string $\sigma$ and any $c\in \N^\N$,
      pass to $\lfloor c\rfloor_k$ and overwrite the first $\abs{\sigma}$ bits
      of the latter string with those of $\sigma$ to produce an instance $d$ of
      $f \restr_{[\sigma]}$. If $H$ is an
      $f$-solution to $d$, then it is also an $\Ind[k]{\calS}$-solution to $d$,
      of course. Since $d$ agrees with $\lfloor c\rfloor_k$ past the first
      $\abs{\sigma}$ bits, $\Phi(H,\abs{\sigma})$ is an
      $\Ind[k]{\calS}$-solution to $\lfloor c\rfloor_k$, hence an $f$-solution
      to the same, and hence also an $f$-solution to the original $c$. 
    
    \item It suffices to show $\Ind[\N]{\calS} \leq_\W \Ind[\N]{\calS}
      \restr_{[\sigma]}$ for every $\sigma \in \N^{<\N}$. Let $n$ be the largest
      entry of $\sigma$. Given any bounded coloring $c$ of $S$, produce a new
      coloring $d \in [\sigma]$ by replacing the first $\abs{\sigma}$ bits of
      $c$ with $\sigma$ and adding $n+1$ to the remaining bits. Then every
      $d$-solution is also a $c$-solution. \qedhere
  \end{enumerate}
\end{proof}

It is known that if $f$ is a fractal and $f \leq_\W \bigsqcup_{i\in \N}
g_i$ for any problems $g_i$, then $f \leq_\W g_i$ for some $i$
\cite[Lemma~5.5]{closedchoice}.
Additionally, if $f$ is a total fractal and $f\leq_\W g\ast \CN$ for any $g$,
then $f\leq_\W g$ \cite[Theorem~2.4]{LRP}. 
Thus $\Ind[\N]{\calS}$ is not in general a total fractal: if it were,
\thref{IndSN} would imply $\Ind[\N]{\calS} \leq_\W \Ind{\calS}$, and hence
$\Ind[\N]{\calS} \leq_\W \Ind[k]{\calS}$ for some $k$.
This fails for $\Eqinf$ by \thref{eqsep1} and for $\Q$ by the discussion in
Section~\ref{sec:prior}, which shows that $\Ind[k+1]{\Q} \not\leq_\W
\Ind[k]{\Q}$.
(It seems unlikely that it could happen for any $\calS$, though we do not have
a proof it does not.) For the same reason, neither $\Ind{\Q}$ nor $\Ind{\Eqinf}$
is a fractal. 
One may compare this situation to that of the $\mathsf{Shuffle}$ family of
problems, which were studied in \cite{PPS23} and whose definition is given in
the next section. As per the discussion in Section~6 of that paper, the versions
of those problems with a fixed number of colors are all total fractals; the
versions where instances have no information on the number of colors are
fractals, but not total fractals; and the versions where instances have a given
but not fixed number of colors are not fractals at all.

The condition of uniform computable categoricity in
\thref{IndSstar,IndSN,SNfractal} is sufficient but not necessary, since
they hold for $\Ind\N \equiv_\W \RT^1_+$ while $(\N,<)$ is not computably
categorical. (It is perhaps worth noting that $\Ind (\N,<)$ is Weihrauch
equivalent to $\Ind (\N,\emptyset)$, the structure consisting of a countably
infinite set in the empty language, and the latter structure \emph{is} u.c.c.)
We do not know a characterization of the structures $\calS$ for
which the above results hold.
A good source of
(counter)examples could be the class of nonscattered (countable) linear orders,
all of which are indivisible as a consequence of the fact that every countable
linear order embeds into $\Q$. The only infinite u.c.c.\ linear order is $\Q$
itself (see e.g.~\cite[Example~III.4]{Mon1}), but the indivisibility problems of
many other linear orders turn out to have exactly the same computational
strength as those of $\Q$.
A linear order $\mathcal{L}$ is called \bdef{strongly $\eta$-like} if there is a
number $n$ such that $\mathcal{L}$ can be built from $\Q$ by replacing every
element by a chain of between one and $n$ elements. (Here $\eta$ is the order
type of $\Q$. See for instance \cite[\S4]{DowneyLO} or \cite[\S10.1]{STT} for a
discussion of this notion.)
Then we have

\begin{prop}\thlabel{nQm} If $\mathcal{M}$ is a strongly $\eta$-like computable
  linear order, and $\mathcal{L} = (m_1+\mathcal{M}+m_2,<)$ for some $m_1,m_2\in
  \N$, then $\Ind[k]{\mathcal{L}} \equiv_\W \Ind[k]{\Q}$. Similarly for
  $\Ind{\mathcal{L}}$ and $\Ind[\N]{\mathcal{L}}$.
\end{prop}
\begin{proof} ~[$\leq_\W$] Given an instance $\ang{\calA, c}$ of
  $\Ind[k]{\mathcal{L}}$, we will produce a subcopy of $\Q$ as follows.
  Let $n$ witness that $\mathcal{M}$ is strongly $\eta$-like.
  Read far enough in the atomic diagram of $\calA$ to find two points $x_0 <
  x_1$ such that there are more than $m_1$ points below $x_0$, more than $m_2$
  above $x_1$, and more than $n$ between $x_0$ and $x_1$. At each subsequent
  stage, for every $x_i < x_j$ found so far, search for an $x_h \in (x_i,x_j)$
  such that $(x_i,x_h)$ and $(x_h,x_j)$ each contain more than $n$ points. This
  criterion ensures that $x_h$ cannot be finitely apart from
  either $x_i$ or $x_j$, so that both
  $(x_i,x_h)$ and $(x_h,x_j)$ contain densely ordered sets. Hence it is always
  possible to find a suitable $x_h$ at every step, and $X = \{x_2,x_3,\dotsc\}$
  is a uniformly $\calA$-computable set which induces a substructure
  $\mathcal{X} \simeq \Q$.
  Then $\ang{\mathcal{X},c}$ is an instance of $\Ind[k]{\Q}$; if $H$ is a
  solution to this instance, then by density of $H$, we can use a greedy
  algorithm to $\calA$-compute a subcopy of $\calA$ within $H$. This subcopy is
  a solution to the original instance $\ang{\calA,c}$ of $\Ind[k]{\mathcal{L}}$.
  (This direction of the reduction does not require $\mathcal{M}$ to have a
  computable copy.)

  [$\geq_\W$] The proof is essentially the same as the above, but in reverse. 
  Let $\ang{\calA,c} \in \dom (\Ind[k]{\Q})$. Build an instance of
  $\Ind[k]{\mathcal{L}}$ by producing a subcopy $\mathcal{B}$ of $\mathcal{L}$
  inside $\calA$ and using $c$ as-is. (This step is the only place where we use
  that $\mathcal{L}$ is computable, since the forward functional in this
  reduction does not have oracle access to a presentation of $\mathcal{L}$.) If
  $H$ is an $\Ind[k]{\mathcal{L}}$-solution to $\ang{\mathcal{B}, c}$, then one
  can carry out the procedure described above to find $X$ inside $H$, and $X$
  will be a solution to $\ang{\calA,c}$.\qedhere
\end{proof}

Consequently, the conclusions of \thref{IndSstar,IndSN,SNfractal} hold
for any such $\mathcal{L}$.
Work is currently in progress to study various weakenings of uniform computable
categoricity with the aim of recovering these and other properties for a much
wider class of structures.


\section{The rational numbers}\label{sec:rats}

\subsection{Prior related work}\label{sec:prior}
$\Ind{\Q}$ is exactly the problem form of the reverse-mathematical principle
$\mathsf{ER}^1$ studied by Frittaion and Patey in \cite{FP17}. There the authors
show, among other things, that the implication $\mathsf{ER}^1 \to \RT^1_+$ over
$\mathsf{RCA}_0$ is strict. This has a kind of uniform counterpart in our
\thref{pcetcor}(i), since trivially $\RT^1_+ \leq_\W \Ind{\calS}$ (and $\RT^1_\N
\leq_\W \Ind[\N]{\calS}$) for any indivisible $\calS$ with a computable
presentation, and this is strict if $\calS = \Q$. 
But separations over $\mathsf{RCA}_0$ have no direct bearing in the
present setting: both $\Q$ and $\N$ are computably indivisible and hence both
$\Ind{\Q}$ and $\RT^1_+$ hold in $\omega$-models of $\mathsf{RCA}_0$, indeed are
computably equivalent to the identity map on Baire space. To obtain meaningful
distinctions between them one must pass to a stronger reducibility, and so we
focus on the uniform content of $\Ind{\Q}$ via its Weihrauch degree.

The Weihrauch degrees of a family of problems related to $\Ind{\Q}$ have
recently been studied by Pauly, Pradic, and Sold\`a
\cite{PPS23}. In their terminology, if $c$ is any coloring of $\Q$, then an open
interval $I\subseteq\Q$ is a \bdef{$c$-shuffle} if for every color occurring in
$I$, the set of points of that color is dense in $I$. They then investigate
several corresponding families of problems:
\begin{itemize}
  \item $\mathsf{Shuffle}$ has instances $(k,c)$, where $k\in\N$ and $c\colon
    \Q\to k$. If $I$ is a (code for a) rational interval and $C\subseteq k$,
    then $(I,C)\in \mathsf{Shuffle}(k,c)$ iff $I$ is a $c$-shuffle with exactly
    the colors of \nolinebreak $C$.
  \item $\mathsf{iShuffle}$ is the weakening of $\mathsf{Shuffle}$ which, for an
    instance $(k,c)$, returns an interval $I$ such that there exists $C\subseteq
    k$ with $(I,C)\in \mathsf{Shuffle}(k,c)$.
  \item $\mathsf{cShuffle}$ is the weakening of $\mathsf{Shuffle}$ which, for an
    instance $(k,c)$, returns a $C\subseteq k$ such that there is an $I$ with
    $(I,C)\in \mathsf{Shuffle}(k,c)$.
  \item $(\eta)^1_{<\infty}$, which had previously been studied in \cite{FP17}
    as a reverse-mathematical principle, is a different weakening of
    $\mathsf{Shuffle}$ which, for an instance $(k,c)$, returns an interval $I$
    and a single color $n<k$ such that the points of $c$-color $n$ are dense in
    $I$.
\end{itemize}

$\mathsf{Shuffle}_k$, $\mathsf{iShuffle}_k$, and $\mathsf{cShuffle}_k$ are
defined to be the restrictions of the first three problems to instances of the form
$(k,c)$. One can also define $\mathsf{Shuffle}_\N$,
$\mathsf{iShuffle}_\N$, and $\mathsf{cShuffle}_\N$ in a natural way.
The authors show that $\mathsf{cShuffle} \equiv_\W (\mathsf{LPO}')^\ast$,
$\mathsf{iShuffle} \equiv_\W \TCN^\ast \equiv_\W (\eta)^1_{< \infty} \equiv_\W
\mathsf{i}(\eta)^1_{<\infty}$ where the latter problem is the version of
$(\eta)^1_{<\infty}$ which only returns an interval $I$, and $\mathsf{Shuffle}
\equiv_\W (\mathsf{LPO}')^\ast \times \TCN^\ast$. 
More specifically, they show $\mathsf{iShuffle}_k \equiv_\W \TCN^{k-1}$, the
$(k-1)$-fold parallel product of $\TCN$. 
(They also show that
$\mathsf{cShuffle}_k \leq_\W (\mathsf{LPO}')^{2^k-2}$. The reverse direction of
the equivalence between $\mathsf{cShuffle}$ and $(\mathsf{LPO}')^\ast$ is
established by showing that $\mathsf{LPO}' \leq_\W \mathsf{cShuffle}$ and
$\mathsf{cShuffle} \equiv_\W \mathsf{cShuffle}^\ast$. The precise relationship
between number of colors and number of parallel instances of $\mathsf{LPO}'$ is
left open.)
As part of their investigation, they establish further that
$\mathsf{iShuffle}_\N \equiv_\W \mathsf{iShuffle} \ast \CN$ (Proposition~63),
and also that $\cRT^1_{k+1} \leq_\W \TCN^m$ if and only if $k\leq m$
(Theorem~10). Since $\RT^1_\N$ is a fractal, it follows that $\mathsf{iShuffle}$
cannot solve $\RT^1_\N$: if it could, then we would have $\RT^1_\N \equiv_\W
\cRT^1_\N \leq_\W \TCN^k$ for some $k$, a contradiction.

It is clear that $\Ind[k]{\Q} \leq_\W \mathsf{iShuffle}_k$, $\Ind{\Q} \leq_\W
\mathsf{iShuffle}$, and $\Ind[\N]{\Q} \leq_\W \mathsf{iShuffle}_\N$, since for
any color $i$ found in a shuffle, the set of points of that color in the shuffle
is isomorphic to $\Q$. 
As $\mathsf{iShuffle}_k$ cannot solve $\cRT^1_{k+1}$, it follows that
$\Ind[k+1]{\Q} \not\leq_\W \Ind[k]{\Q}$.
Similarly, the fact that $\RT^1_\N$ is below $\Ind[\N]{\Q}$ but not
$\mathsf{iShuffle}$ implies that $\Ind[\N]{\Q} \not\leq_\W \mathsf{iShuffle}$
and in particular $\Ind[\N]{\Q} \not\leq_\W \Ind{\Q}$. 

After a draft of this article appeared, the author was made aware of
contemporaneous work by Dzhafarov, Solomon, and Valenti \cite{DSV} concerning
the \bdef{tree pigeonhole principle} $\mathsf{TT}^1_+$. This states that given
any coloring of $2^{<\N}$ with bounded range, there is an infinite monochromatic
subset $H$ isomorphic to $2^{<\N}$ as a partial order---in other words, that
$(2^{<\N},\prec)$ is indivisible, where $\prec$ is the prefix relation. Problems
$\mathsf{TT}^1_k$ and $\mathsf{TT}^1_\N$ are defined analogously. An earlier
version of this paper claimed that $\mathsf{TT}^1_k \equiv_\W \Ind[k]{\Q}$,
$\mathsf{TT}^1_+ \equiv_\W \Ind{\Q}$, and $\mathsf{TT}^1_\N \equiv_\W
\Ind[\N]{\Q}$. While the problems do indeed appear to be very closely related,
the proof which was sketched of their claimed equivalence was erroneous, and it
seems likely that $\Ind{\Q}$ is strictly Weihrauch reducible to
$\mathsf{TT}^1_+$.
However, the question of their separation remains open and we leave its
resolution to future work. 
We would like to thank Reed Solomon for bringing the aforementioned error to
light.

Among the results of \cite{DSV} are several corollaries whose conclusions are
similar to ours, but they are obtained by entirely different methods, and there
is otherwise no overlap with the present work.
On the other hand, Pauly \cite{P24} has independently proven a stronger
version of \thref{QperpCN} for $\mathsf{TT}^1_+$, namely that
$\mathsf{C}_{k+1} \not\leq_\W \mathsf{TT}^1_k$ where $\mathsf{C}_k$ outputs an
element of a nonempty subset of $k = \{0,\dotsc, k-1\}$ given an enumeration of
its complement. 
His paper also contains an improved version of our \thref{pcet} below for
$\mathsf{TT}^1_+$, recovering it as part of a more general investigation. We
present our proofs anyway in hope of the ideas therein being useful elsewhere,
perhaps in non-reductions involving other indivisibility problems.

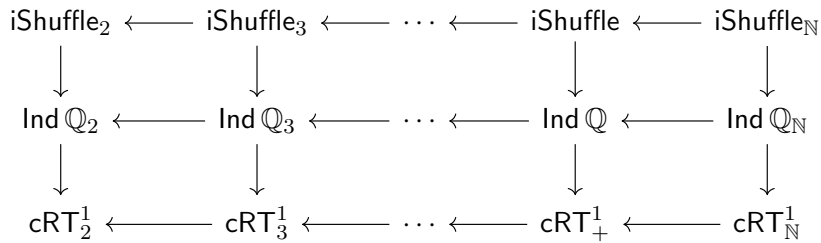
\begin{figure}
  \centering
  \begin{tikzcd}
    \mathsf{iShuffle}_2 \arrow[d] & \mathsf{iShuffle}_3 \arrow[l] \arrow[d] &
    \cdots \arrow[l] &
    \mathsf{iShuffle} \arrow[l] \arrow[d] & \mathsf{iShuffle}_\N \arrow[l]
    \arrow[d]\\
    \Ind[2]{\Q} \arrow[d] & \Ind[3]{\Q} \arrow[l] \arrow[d] &
    \cdots \arrow[l] & \Ind{\Q} \arrow[l]
    \arrow[d] & \Ind[\N]{\Q} \arrow[l] \arrow[d]\\
    \cRT^1_2 & \cRT^1_3 \arrow[l] & \cdots \arrow[l] &
    \cRT^1_+ \arrow[l] & \cRT^1_\N \arrow[l]
  \end{tikzcd}
  \caption{Reductions between versions of $\cRT^1$, $\Ind{\Q}$, and
      $\mathsf{iShuffle} \equiv_\W \TCN^\ast$. An arrow
    from $Q$ to $P$ signifies that $P<_\W Q$. No other reductions hold other than
  those implied by transitivity.}
  \label{Qreds}
\end{figure}

\subsection{The Weihrauch degree of \texorpdfstring{$\Ind{\Q}$}{Ind Q}}

The main result of this section is that the reduction of $\Ind[k]{\Q}$ to
$\mathsf{iShuffle}_k$ is strict (\thref{Qstrict}), which follows from
\thref{QperpCN} below as $\CN$ is Weihrauch reducible to $\TCN \equiv_\W
\mathsf{iShuffle}_2$.
The relationships between $\Ind[k]{\Q}$, $\TCN^k$, and $\cRT^1_k$ are summarized
in Figure~\ref{Qreds}. 
The fact that $\mathsf{iShuffle}_2 \not\leq_\W \Ind[\N]{\Q}$ is another
consequence of the theorem (\thref{noTCN}), and we have $\Ind[2]{\Q} \not\leq_\W
\cRT^1_\N$ by \thref{pcetcor} of \thref{pcet} below.
All other (non)reductions in the figure were justified in the previous section.

We will identify without comment natural numbers with the rationals they encode
via a fixed computable presentation of $\Q$.

\begin{thm}\thlabel{QperpCN} $\CN \not\leq_\W \Ind{\Q}$.
\end{thm}

\begin{proof} 
  It is known that $\CN$ is a fractal, so if $\CN \leq_\W \Ind{\Q} = \bigsqcup_k
  \Ind[k]{\Q}$ then $\CN \leq_\W \Ind[k]{\Q}$ for some $k$. It suffices to show
  this cannot occur. Fix $k$ and suppose for the sake of contradiction that $\CN
  \leq_\W \Ind[k]{\Q}$ via $\Delta$ and $\Psi$.
  
  Informally, the strategy is to wait for $\Psi$ to output a number using some
  finite monochromatic set $F$ as an oracle, and then try to diagonalize by
  enumerating that number.
  If we can extend from there to a $\CN$-instance for which $F$ is still part of
  a valid $\Ind[k]{\Q}$-solution, we are done. If not, then no matter how we
  extend, there will be an $F$-interval---an interval in the partition
  of $\Q$ induced by $F$---in which red points are scattered, if $F$ was red.
  The particular $F$-interval where this is the case depends on the
  $\CN$-instance, so we account for all possibilities,
  repeating the same procedure as before in every $F$-interval,
  now looking for $\Psi$ to converge on a blue set in each and diagonalizing
  again by enumerating all of its outputs into our instance of $\CN$. 
  This completes the proof if $k=2$ since blue points must be dense in any
  interval where red points are scattered.
  If there are more than two colors, and if again none of the blue sets found
  can possibly form part of a valid $\Ind[k]{\Q}$-solution no matter how we
  continue to extend, then there will
  always be some subinterval induced by the blue sets in which blue points are
  scattered, and in particular there will be one in which \emph{both} red and
  blue points are scattered. This is illustrated in Figure~\ref{cnfig}.
  Since every interval has a subinterval in which some color is dense,
  if the procedure is iterated, then eventually $\Delta$ cannot force any more
  colors to be scattered and so the final diagonalization attempt succeeds.

  Now we proceed to the formal proof.
  We regard $\Psi$ as computing a $\CN$-solution $n$ if $\Psi^{ f\oplus
  H}(0)\cvg =n$, for an instance $f$ of $\CN$ and $\Ind[k]{\Q}$-solution $H$ of
  $\Delta^f$. A solution $n$ is valid only if $n\notin \ran f$.
  Begin by finding a string $\sigma_0 \in \N^{<\N}$ and a finite set $S_\lambda
  \subset \Q$ such that $\Delta^{\sigma_0}(S_\lambda) \cvg = 0$ (without loss of
  generality), and such that $\Psi^{\sigma_0 \oplus S_\lambda}(0) \cvg =
  m_\lambda$ for some $m_\lambda$. Here $\lambda$ denotes the empty string. 
  One can obtain $\sigma_0$ and $S_\lambda$ as initial segments of any $f\in
  \dom\CN$ and any $H$ solving $\Delta^f$. If $m_\lambda \in \ran \sigma_0$, the
  reduction fails already, so assume otherwise and let $\sigma_0' =
  \sigma_0^\smallfrown m_\lambda$ be our first (or to be consistent later,
  zeroth) attempt at diagonalizing.

  Suppose first there is a $g\in [\sigma_0'] \cap \dom \CN$ such
  that $S_\lambda$ is extendible to a solution $K$ of $\Delta^g$. If so, then
  $\Psi^{g\oplus K}(0)\cvg = m_\lambda$ but $m_\lambda \in \ran g$, defeating
  the reduction in this case. If not, write the elements of $S_\lambda$ in
  increasing order as $x^\lambda_1 < \dotsc < x^\lambda_{\ell(\lambda)}$ where
  $\ell(\lambda) = \abs{ S_\lambda}$; 
  let $x^\lambda_0 = -\infty$ and $x^\lambda_{\ell(\lambda)+1} = +\infty$; 
  let $X_\lambda = S_\lambda \cup \{x^\lambda_0, x^\lambda_{\ell(\lambda)+1}\}$;
  and let $I^\lambda_i$ be the interval $(x^\lambda_i, x^\lambda_{i+1})$.
  If $S_\lambda$ is not extendible to a solution of $\Delta^g$ for
  any $g\in [\sigma_0'] \cap \dom \CN$, then for each $g$ there is at least one
  $i \geq 0$ such that the set of color-0 points in $I^\lambda_i$ with respect
  to $\Delta^g$ is scattered. 
  Then the complementary set of points of $\Delta^g$-colors $1,2,\dotsc,k-1$ is
  dense in that interval.
  We are then guaranteed to find monochromatic sets of at least one color $1,2,
  \dotsc, k-1$ in each such interval which witness convergence of $\Psi$ on some
  extension of $\sigma_0'$, and can repeat the diagonalization strategy on each,
  using a procedure described in the next paragraph.
  This will produce a tree structure of nested intervals labeled by finite
  strings $\beta$ as displayed in Figure~\ref{cnfig}.

  Let $L>0$ and suppose $\sigma_{L-1}'$ has already been found, and that the sets
  $X_\beta$ and intervals $I^\beta_i$ have been defined for all $\abs{ \beta}
  = L-1$ and $0\leq i \leq \ell(\beta)$.
  At this stage, $\sigma_{L-1}'$ represents the $(L-1)$st attempt---starting
  with the $0$th---to diagonalize
  against the reduction, having already enumerated the outputs of $\Psi$ on each
  of the sets $S_\delta$ for all $\abs{\delta} \leq L-1$ for which such a set
  exists.
  Each $S_\delta$ has color $\abs{\delta}$ and was found inside the interval
  $I^\gamma_i$ where $\delta = \gamma^\smallfrown i$, or simply inside
  $(-\infty,\infty)$ if $\delta = \lambda$ (corresponding to the ``base case''
  $L=0$ which was described in the last paragraph).
  We claim that there is a $\sigma_L \in [\sigma_{L-1}']$ which, intuitively
  speaking, sees $\Psi$ converge on a set of color $L$ in as many intervals
  $I^\beta_i$ as possible, meaning no extension of $\sigma_L$ does so in
  any additional interval.
  To be precise,
  for all $\rho\in [\sigma_L]$, if any $i\geq0$ and $\beta$ with
  $\abs{\beta}=L-1$ are such that there exists a finite set $S \subset
  I^\beta_i$ which has $\Delta^\rho$-color $L$ and such that $\Psi^{\rho \oplus
  S}(0)\cvg$, then the same is true of $\sigma_L$ in place of $\rho$, for the
  same $i$ (and possibly a different $S$).
  To prove the claim, consider the function $a \colon [\sigma_{L-1}'] \to \N$
  which sends a string $\tau$ to the number of distinct intervals $I^\beta_i$,
  for $\abs{\beta} = L-1$ and $0\leq i\leq \ell(\beta)$, such that there exists
  $S_{\beta^\smallfrown i}\subset I^\beta_i$ witnessing $\Delta^\tau
  (S_{\beta^\smallfrown i}) \cvg = L$ and $\Psi^{\tau \oplus
  S_{\beta^\smallfrown i}}(0)\cvg$. Since there are only finitely many intervals
  $I^\beta_i$ at hand, $a$ is bounded, and any $\sigma_L \in [\sigma_{L-1}']$
  witnessing its maximum value has the stated property.

  \begin{sidewaysfigure}
     \centering
\vskip15cm
    \begin{tikzpicture}[scale=0.4]
      \draw (-25,2.5) node{$L=0$}
            (-25,-2.5) node{$L=1$}
            (-25,-7.5) node{$L=2$}
            (-25,-12.5) node{$\vdots$};
      \draw [dashed]
            (-21,-0.5) -- (-21,-15)
            (0,0) -- (0,-15)
            (25.5,-0.5) -- (25.5,-15);
      \draw [dashed]
            (-10,-5.3) -- (-10,-15.2)
            (5,-5.3) -- (5,-15.2)
            (15,-5.3) -- (15,-15.2);
      \draw [dashed]
            (-18,-10) -- (-18,-15)
            (-14,-10) -- (-14,-15)
            (2,-10) -- (2,-15)
            (7,-10) -- (7,-15)
            (12,-10) -- (12,-15)
            (21,-10) -- (21,-15);

      \draw (-10,0) node{\tiny{$I^\lambda_0$}}
            (12.5,0) node{\tiny{$I^\lambda_1$}}
            (-15,-5) node{\tiny{$I^0_0$}}
            (-5,-5) node{\tiny{$I^0_1$}}
            (2.5,-5) node{\tiny{$I^1_0$}}
            (10,-5) node{\tiny{$I^1_1$}}
            (20,-5) node{\tiny{$I^1_2$}}
            (-19.5,-10) node{\tiny{$I^{00}_0$}}
            (-16,-10) node{\tiny{$I^{00}_1$}}
            (-12,-10) node{\tiny{$I^{00}_2$}}
            (-5,-10) node{\tiny{$I^{01}_0$}}
            (1,-10) node{\tiny{$I^{10}_0$}}
            (3.5,-10) node{\tiny{$I^{10}_1$}}
            (6,-10) node{\tiny{$I^{11}_0$}}
            (9.5,-10) node{\tiny{$I^{11}_1$}}
            (13.5,-10) node{\tiny{$I^{11}_2$}}
            (18,-10) node{\tiny{$I^{12}_0$}}
            (23,-10) node{\tiny{$I^{12}_1$}};

      \tikzset{every path/.style={line width=0.05mm}}
      \draw [->] (-11,0) -- (-20.5,0); 
      \draw [->] (-9,0) -- (-1,0);
      \draw [->] (11.5,0) -- (1,0); 
      \draw [->] (13.5,0) -- (25,0);
      \draw [->] (-16,-5) -- (-20.5,-5); 
      \draw [->] (-14,-5) -- (-11,-5);
      \draw [->] (-6,-5) -- (-9,-5); 
      \draw [->] (-4,-5) -- (-1,-5);
      \draw [->] (1.5,-5) -- (0.5,-5); 
      \draw [->] (3.5,-5) -- (4.5,-5);
      \draw [->] (9,-5) -- (6,-5); 
      \draw [->] (11,-5) -- (14,-5);
      \draw [->] (19,-5) -- (16,-5); 
      \draw [->] (21,-5) -- (25,-5);
      \draw [->] (-20,-10) -- (-20.5,-10); 
      \draw [->] (-19,-10) -- (-18.5,-10);
      \draw [->] (-16.5,-10) -- (-17.5,-10); 
      \draw [->] (-15.5,-10) -- (-14.5,-10);
      \draw [->] (-12.5,-10) -- (-13.5,-10); 
      \draw [->] (-11.5,-10) -- (-10.5,-10);
      \draw [->] (-6,-10) -- (-9,-10); 
      \draw [->] (-4,-10) -- (-1,-10);
      \draw [->] (0.5,-10) -- (0.25,-10); 
      \draw [->] (1.5,-10) -- (1.75,-10);
      \draw [->] (3,-10) -- (2.3,-10); 
      \draw [->] (4,-10) -- (4.7,-10);
      \draw [->] (5.5,-10) -- (5.2,-10); 
      \draw [->] (6.4,-10) -- (6.75,-10);
      \draw [->] (9,-10) -- (7.5,-10); 
      \draw [->] (10,-10) -- (11.5,-10);
      \draw [->] (13,-10) -- (12.5,-10); 
      \draw [->] (14,-10) -- (14.5,-10);
      \draw [->] (17,-10) -- (15.5,-10); 
      \draw [->] (19,-10) -- (20.5,-10);
      \draw [->] (22.5,-10) -- (21.5,-10); 
      \draw [->] (23.5,-10) -- (25,-10);

      \tikzset{every path/.style={line width=0.4pt}}

      \draw [red!70] (0,2.5) node{$S_\lambda$};
      \filldraw [red!70] (0,0) circle (5pt) node [above] {\tiny{$x_1^\lambda$}};
      \draw (-21,0) circle (5pt) node [above,red!70] {\tiny{$x_0^\lambda = -\infty$}}
        (25.5,0) circle (5pt) node [above,red!70] {\tiny{$x_2^\lambda = \infty$}};
      \draw (-21,-5) circle (2pt)
          (-21,-10) circle (2pt)
          (25.5,-5) circle (2pt)
          (25.5,-10) circle (2pt);
      \filldraw (0,-5) circle (2pt)
          (0,-10) circle (2pt);

      \draw [blue!60] (-10,-2.5) node{$S_0$}
                      (12.5,-2.5) node{$S_1$};
      \filldraw [blue!60] (-10,-5) circle (5pt) node [above] {\tiny{$x_1^0$}}
        (5,-5) circle (5pt) node [above] {\tiny{$x_1^1$}}
        (15,-5) circle (5pt) node [above] {\tiny{$x_2^1$}};
      \filldraw (-10,-10) circle (2pt)
        (5,-10) circle (2pt)
        (15,-10) circle (2pt);

      \draw [orange!70] (-15,-7.5) node{$S_{00}$} 
        (-5,-7.5) node{$S_{01} = \emptyset$} 
        (2.5,-7.5) node{$S_{10}$}
        (10,-7.5) node{$S_{11}$}
        (20,-7.5) node{$S_{12}$};
      \filldraw [orange!70] (-18,-10) circle (5pt) node[above] {\tiny{$x_1^{00}$}}
        (-14,-10) circle (5pt) node[above] {\tiny{$x_2^{00}$}}
        (2,-10) circle (5pt) node[above] {\tiny{$x_1^{10}$}}
        (7,-10) circle (5pt) node[above] {\tiny{$x_1^{11}$}}
        (12,-10) circle (5pt) node[above] {\tiny{$x_2^{11}$}}
        (21,-10) circle (5pt) node[above] {\tiny{$x_1^{12}$}};
    \end{tikzpicture}
  \caption{An illustration of an example diagonalization, with
      rational points vertically displaced to show which are added at each stage $L$.
      Stages are labeled by colors, with $L=0$ being red, $L=1$ blue, and $L=2$
      orange here.
      The point $x^\lambda_1$ is carried through to become the largest element
      of $X_0$ and $X_{01}$ and the smallest element of $X_1$ and $X_{10}$.
      One carries through $x^0_1$, $x^1_1$, and $x^1_2$ in a similar fashion, as
      well as $x^\lambda_0 = -\infty$ and $x^\lambda_2 = +\infty$. The latter
      two points are formally added merely to avoid having to treat the
      outermost intervals as special cases.
      In this example, no orange subset of $I^0_1$ could be found witnessing
      convergence of $\Psi$ on any extension of $\sigma_1$, so we set $X_{01} =
      \{ x^0_1, x^\lambda_1\} = \{ x^0_1, x^0_2 \}$.
      At the end of the procedure, at least one of $I^\lambda_0,I^\lambda_1$
      will have red points scattered; at least one of $I^0_0,I^0_1$ and one of
      $I^1_0,I^1_1,I^1_2$ will have blue points scattered; $I^{01}_0$ will have
      orange points scattered, as will at least one of the $I^{00}_i$s, one of
      the $I^{10}_i$s, one of the $I^{11}_i$s, and one of the $I^{12}_i$s; and
      so on.
  }\label{cnfig}
\end{sidewaysfigure}
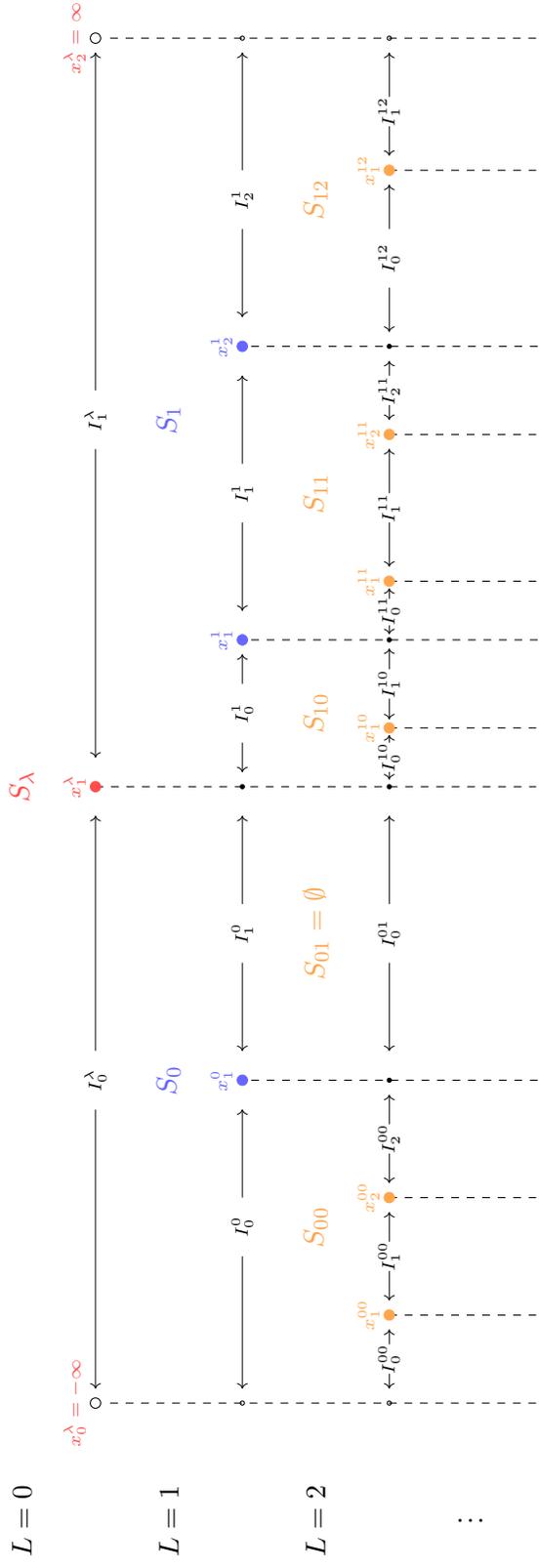

  Continue by picking any such $\sigma_L$, and fix a choice of $S_\alpha \subset
  I^\beta_i$ as
  above for any $\alpha = \beta^\smallfrown i$ for which a set $S_\alpha$ could
  be found. Extend $\sigma_L$ to $\sigma_L'$ by enumerating the output
  $m_\alpha$ of $\Psi^{\sigma_L \oplus S_\alpha}(0)$ for each such $\alpha$.
  Write $S_\alpha = \{ x^\alpha_1 < x^\alpha_2 < \cdots <
  x^\alpha_{\ell(\alpha)} \}$, let $x^\alpha_0 = x^\beta_i < x^\alpha_1$ and
  $x^\alpha_{\ell(\alpha)+1} = x^\beta_{i+1} > x^\alpha_{\ell(\alpha)}$, and let
  $X_\alpha = S_\alpha \cup \{ x^\alpha_0, x^\alpha_{\ell(\alpha)+1} \}$. 
  If for some $\alpha = \beta^\smallfrown i$ of length $L$ no $S_\alpha$ could
  be found, simply let $\ell(\alpha) = 0$ and $X_\alpha = \{ x^\alpha_0,
  x^\alpha_1 \} = \{ x^\beta_i, x^\beta_{i+1} \}$. 
  Then let $I^\alpha_i = (x^\alpha_i, x^\alpha_{i+1})$ for all $\alpha$ of
  length $L$ and all $i=0, \dotsc, \ell(\alpha)$. 
  Thus if no set $S_\alpha$ as above could be found for $\alpha =
  \beta^\smallfrown i$, then we have $I^\alpha_0 = I^\beta_i$, so at the next
  stage one searches for an $S_{\alpha^\smallfrown 0}$ of color $L+1$ in the
  same interval rather than splitting it further.

  If there is a $g\in [\sigma_L'] \cap \dom\CN$ and an $\alpha$ of length $L$
  such that $S_\alpha$ is extendible to a color-$L$ solution of $\Delta^g$, then
  this $g$ witnesses the failure of the reduction since for some solution
  $U\supset S_\alpha$ of $\Delta^g$, we have $\Psi^{g\oplus U}(0)\cvg =
  m_\alpha$ while $m_\alpha \in \ran g$. Suppose, then, that no such $g$ exists.
  This means that for all $g\in [\sigma_L'] \cap \dom \CN$, for all $\beta$ with
  $\abs{\beta} = L-1$, there is an $i\geq0$ such that the set of color-$L$
  points in $I^\beta_i$ is scattered with respect to $\Delta^g$. If $L<k-1$,
  then continue inductively as above by finding sets $S_\xi$ for $\abs{\xi} =
  L+1$ and a string $\sigma_{L+1} \supset \sigma_L'$. But if $L=k-1$, then
  we claim that in fact for any $g\in [\sigma_L'] \cap \dom \CN$, there must be
  an $S_\alpha$ with $\abs{\alpha}=L$ that is extendible to a color-$(k-1)$
  solution of $\Delta^g$. 
  To see why, recall that by assumption at this
  stage, for any such $g$, in particular there is an $i_0$ such that color-0
  points are scattered in $I^\lambda_{i_0}$. In turn, there is an $i_1$ such
  that color-1 points are also scattered in $I^{i_0}_{i_1}$; an $i_2$ such that
  color-2 points are also scattered in $I^{i_0i_1}_{i_2}$; and so on, so that
  there is ultimately an $i_{k-2}$ such that the points of colors $0,1,\dotsc,
  k-2$ are all scattered in the interval $I^\beta_{i_{k-2}}$, where $\beta =
  i_0i_1 \cdots i_{k-3}$ (or $\beta=\lambda$ if $k=2$). It follows that the set
  of points of color $k-1$ is dense in $I^\beta_{i_{k-2}}$. In particular, the
  set of color-$(k-1)$ points in $I^\beta_{i_{k-2}}$ is a solution of
  $\Delta^g$, so if $\alpha = i_0i_1\dotsc i_{k-2}$, then during the procedure
  we must have been able to find a set $S_\alpha$ of color $k-1$ with
  $\Psi^{\sigma_L \oplus S_\alpha}(0)$ converging.
  Hence for this $g$ and $\alpha$, there is a solution $U\supset S_\alpha$ of
  $\Delta^g$ with $\Psi^{g\oplus U}(0)\cvg = m_\alpha$ while $m_\alpha\in \ran
  g$, defeating the reduction and completing the proof. \qedhere
\end{proof}

\begin{cor}\thlabel{Qstrict} For all $k\geq2$, the Weihrauch reducibility of
  $\Ind[k]{\Q}$ to $\mathsf{iShuffle}_k$ is strict.\qed
\end{cor}

\begin{cor}\thlabel{noTCN} 
  $\mathsf{iShuffle}_2 \not\leq_\W \Ind[\N]{\Q}$, and hence in
  particular $\mathsf{iShuffle}_\N \not\leq_\W \Ind[\N]{\Q}$.
\end{cor}
\begin{proof} Recall that $\mathsf{iShuffle}_2 \equiv_\W \TCN$. The latter
  problem is known to be a total fractal.
  We have by \thref{IndSN} that $\Ind[\N]{\Q} \equiv_\W \Ind{\Q} \ast \CN$,
  and it follows that if $\TCN \leq_\W \Ind[\N]{\Q}$ then $\TCN \leq_\W
  \Ind{\Q}$, contradicting \thref{QperpCN}.
\end{proof}

\thref{QperpCN} and its corollaries show that $\Ind{\Q}$ is uniformly rather weak,
but in a sense it is not too easy to solve either, in that it cannot be solved
by any problem having a c.e.\ approximation in the following sense. 
Let $\Phi$ be a universal functional.

\begin{defi}\thlabel{pcetdef} We say a problem $P$ is \bdef{pointwise c.e.\
  traceable}, or p.c.e.t.\ for short, if there is an index $i$
  such that for all $p\in \dom P$, $\Phi_{\ang{i,p}}(0)$ converges and outputs
  some $e$ such that $W_e^p$ is finite and contains at least one index $j$ with
  $\Phi_{\ang{j,p}} \in P(p)$.
\end{defi}
Intuitively, this means there is a uniform procedure to enumerate a finite list of
potential solutions to any instance of $P$, at least one of which will turn
out to be correct.
Examples of pointwise c.e.\ traceable problems include $\RT^1_\N$ and any
problem with finite computable codomain.

\begin{thm}\thlabel{pcet} If $P$ is any pointwise c.e.\ traceable problem,
  then $\Ind[2]{\Q} \not\leq_\W P$.
\end{thm}

Hence none of $\Ind[k]{\Q}$, $\Ind{\Q}$, and $\Ind[\N]{\Q}$ are reducible to $P$
either.
Although this result was found independently in the course of the present work,
the core idea of its proof was outlined in \cite[\S5]{FP17}, where the authors
call it a ``disjoint extension commitment'' of $\mathsf{ER}^1$ (or $\Ind{\Q}$),
and may have appeared elsewhere. The idea is that if $\Ind[2]{\Q} \leq_\W P$,
then when the return functional $\Psi$ outputs any points $x<y<z$, it commits to
outputting monochromatic densely ordered sets both in $(x,y)$ and in $(y,z)$.
If there are only finitely many possible solutions that need be considered as
oracles for $\Psi$, then once $\Psi$ has output enough points per solution, we
can choose a finite set of pairwise disjoint intervals, one per solution,
and diagonalize separately
in each by filling it in with the opposite color as its endpoints. This works as
long as we know the number of solutions of
the computed $P$-instance will not grow without bound---this is what the p.c.e.t.\
condition is used for.

\begin{proof}[Proof of \thref{pcet}] Suppose $\Ind[2]{\Q}\leq_\W P$ via $\Delta$
  and $\Psi$. We reach a contradiction by building a coloring $c = \lim_s c_s$
  of $\Q$ by finite extension to witness the failure of this reduction,
  beginning with $c(0)=0$.

  Let $i$ be an index associated to $P$ as in \thref{pcetdef}, and let
  $N(s)$ be either $W^{\Delta(c)}_x[s]$ if $\Phi_{\ang{i,\Delta(c)}}(0)[s] \cvg = x$, or
  $\emptyset$ if this computation diverges. 
  For each $j\in N(s)$, let $\Psi^c_j = \Psi(\Phi_{\ang{j,\Delta(c)}} \oplus
  c)$, the computation of the return functional using the $j$th
  $\Delta(c)$-solution furnished by $\Phi$ as an oracle, along with the original
  coloring $c$.
  Since there is always some $j\in \bigcup_s N(s)$ such that
  $\Phi_{\ang{j,\Delta(c)}}$ solves $\Delta(c)$, at least one $\Psi^c_j$ must
  be a solution to $c$ in the end, so it is enough to diagonalize against every
  $\Psi^c_j$ without regard to its eventual correctness.

  Define a recurrence relation $r(n)$ by letting $r(1) = 3$ and $r(n+1) =
  2nr(n)+2$ for $n\geq1$. 
 
  \begin{lem}\thlabel{pcetlem}
    If sets $S_0,\dotsc, S_{k-1}$ each contain at least $r(k)$ points, then
    there are pairwise disjoint rational intervals $I_0, \dotsc, I_{k-1}$ such that the
    endpoints of $I_i$ are in $S_i$ for each $i$.
  \end{lem}
  \begin{proof}
    By induction on $k$. The case $k=1$ is
    immediate, so suppose the claim holds for all $i\leq k$ and that we have
    $k+1>1$ sets $S_0, \dotsc, S_{k-1}, S_k$ each with
    at least $r(k+1)$ points. For purposes of the inductive step, we
    discard all but $2r(k)$ points from the first $k$ sets $S_0, \dotsc, S_{k-1}$. 
    Then the $2kr(k)$ total points in $S_0 \cup \cdots \cup S_{k-1}$
    divide $\Q$ into $2kr(k)+1$ open intervals. If $S_k$ has at least
    $2kr(k)+2$ points, by the pigeonhole principle, at least one of those
    $2kr(k)+1$ open intervals will contain two distinct points $a_k,b_k\in S_k$.
    We in particular choose $a_k$ and $b_k$ to be adjacent among the elements of
    $S_k$, then set $I_k=(a_k,b_k)$. 
    This splits $\bigcup_{i<k} S_i$
    into two groups, the points below $a_k$ and those above $b_k$.
    Using the pigeonhole principle again, for each $i<k$, one of those groups
    contains at least $2r(k)/2 = r(k)$ elements of $S_i$. The
    inductive hypothesis can now be applied to the collection of (at most $k$)
    sets having at least $r(k)$ points below $a_k$, and separately
    to the collection of sets having at least $r(k)$ points above $b_k$,
    to produce a set of disjoint rational intervals as required.
  \end{proof}

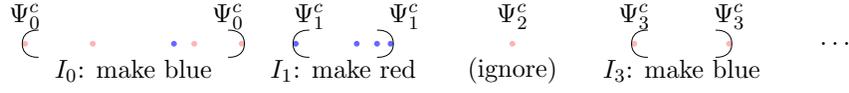
\begin{figure}
\centering
\scalebox{0.9}{
    \begin{tikzpicture}[scale=0.2]
      \filldraw [red!30] (0,0) circle (5pt);
      \draw plot [smooth, tension=2] coordinates {(1,1) (-0.2,0) (1,-1)};
      \filldraw [red!30] (5,0) circle (5pt);
      \filldraw [blue!60] (11,0) circle (5pt);
      \filldraw [red!30] (12.5,0) circle (5pt);
      \filldraw [red!30] (16,0) circle (5pt);
      \draw plot [smooth, tension=2] coordinates {(15,1) (16.2,0) (15,-1)};
      \draw (0,2) node{$\Psi^c_0$};
      \draw (15,2) node{$\Psi^c_0$};
      \draw (8,-2) node{$I_0$: make blue};
      \filldraw [blue!60] (20,0) circle (5pt);
      \draw plot [smooth, tension=2] coordinates {(21,1) (19.8,0) (21,-1)};
      \filldraw [blue!60] (24.5,0) circle (5pt);
      \filldraw [blue!60] (26,0) circle (5pt);
      \filldraw [blue!60] (27,0) circle (5pt);
      \draw plot [smooth, tension=2] coordinates {(26,1) (27.2,0) (26,-1)};
      \draw (21,2) node{$\Psi^c_1$};
      \draw (28,2) node{$\Psi^c_1$};
      \draw (23.5,-2) node{$I_1$: make red};
      \filldraw [red!30] (36,0) circle (5pt);
      \draw (36,2) node{$\Psi^c_2$};
      \draw (36,-2) node{(ignore)};
      \filldraw [red!30] (45,0) circle (5pt);
      \filldraw [red!30] (52,0) circle (5pt);
      \draw plot [smooth, tension=2] coordinates {(46,1) (44.8,0) (46,-1)};
      \draw plot [smooth, tension=2] coordinates {(51,1) (52.2,0) (51,-1)};
      \draw (45,2) node{$\Psi^c_3$};
      \draw (52,2) node{$\Psi^c_3$};
      \draw (48.5,-2) node{$I_3$: make blue};

      \draw (60,0) node{$\cdots$};
    \end{tikzpicture}
  }
  \caption{An example diagonalization against functionals
    $\Psi^c_0,\Psi^c_1,\dotsc$. The points shown have already been colored red
    or blue by $c$ at this stage, and enumerated by $\Psi^c_i$ for some $i \in
    N(s)$.
    $\Psi_0^c$, $\Psi_1^c$, and $\Psi_3^c$ have enumerated enough points so that
    we can find three disjoint intervals $I_0,I_1,I_3$ whose endpoints are
    respectively among the domains of $\Psi_0^c,\Psi_1^c,\Psi_3^c$ up to this
    point. Since the endpoints of $I_0$ are red, we plan to make all future
    points in $I_0$ blue, and similarly for $I_1$ and $I_3$. At this stage,
    $\Psi_2^c$ has not enumerated enough points to find an $I_2$ disjoint from
  $I_0\cup I_1\cup I_3$, so we ignore it for now.} 
  \label{peffigure}
\end{figure}

  Now we proceed to the main argument, which is by finite injury
  (without a priority order). Let $A(s)$ be the set of $i\in N(s)$ such
  that $\Psi_i^c[s]$ has output at least one point. Let $V(s)$ be the set of
  $i\in A(s)$ such that $\Psi_i^c[s]$ has output at least $r(\abs{A(s)})$
  points. At stage $s+1$, for each $i\in V(s)$, if $\Psi_i$ does not have a
  follower, assign it a follower which is a (code for a) rational interval $I_i$
  as furnished by \thref{pcetlem}, applied to the sets $S_i
  = \Psi^c_i[s]$ for $i \in V(s)$.
  Thus the endpoints of $I_i$ were
  enumerated by $\Psi_i^c$ by stage $s$, and $I_i\cap I_j = \emptyset$ if $i\neq
  j\in V(s)$. If there is any $j\in A(s+1)\setminus A(s)$, cancel the followers of
  all $\Psi_i$ which have one at stage $s+1$. Then if $s+1$ is an element of
  $I_i$ for some $i\in V(s)$ and $a$ is an endpoint of $I_i$, color
  $c_{s+1}(s+1) = 1-c_s(a)$. Otherwise, if $s+1$ is not an element of any
  interval $I_i$, color $c_{s+1}(s+1)=0$. (See Figure~\ref{peffigure} for an
  example of how this might look at a particular stage.) 
 
  Because $A(s)$ eventually stabilizes, each $\Psi_i$ will have a follower
  canceled at most finitely many times, and so the set of intervals $I_i$ also
  stabilizes. Any $\Psi_i$ which outputs infinitely many points will have $i\in
  V(s)$ for large enough $s$, so for each such $i$, there are two points output
  by $\Psi_i^c$ between which only finitely many are of the same color.
  Therefore $\Psi_i^c$ does not compute a solution to $c$. \qedhere
\end{proof}

Note that the construction above can be done computably.

\begin{cor}\thlabel{pcetcor} \hfill
  \begin{enumerate}[label=(\roman*)]
  \item $\Ind[2]{\Q} \not\leq_\W \cRT^n_\N$ for all $n$.
  \item $\Ind[2]{\Q} \not\leq_\W (\mathsf{LPO}^{(n)})^\ast$ for all $n$.
  \item $\Ind[2]{\Q} \not\leq_\W \mathsf{cShuffle}$.
\end{enumerate}
\end{cor}
\begin{proof} \hfill
  \begin{enumerate}[label=(\roman*)]
    \item $\cRT^n_\N$ is p.c.e.t.\ since the set of solutions to $c$ is
      contained in $\ran c$. 
    \item If there are $k$ instances of $\mathsf{LPO}^{(n)}$ given in parallel, at
      most $2^k$ distinct solutions are possible, and they can each be encoded as
      a single natural number. So in particular $(\mathsf{LPO}^{(n)})^\ast$ is
      p.c.e.t. 
    \item This follows from (ii) as $\mathsf{cShuffle} \equiv_\W
      (\mathsf{LPO}')^\ast$ \cite{PPS23}. \qedhere
  \end{enumerate}
\end{proof}

Before continuing to the next section, we pause to make some observations about
p.c.e.t.\ problems. Pointwise c.e.\ traceability can be viewed as a
generalization of c.e.\ traceability, which is the special case where $P$ is a
single-valued function from $\N$ to $\N$ (see for instance \cite[\S 11.4]{Soare}
for the definition). It is also related to the notion of a pointwise finite
problem, defined in \cite{GPV21}, which is a problem such that every instance
has finitely many solutions. However, a p.c.e.t.\ problem is not necessarily
pointwise finite, one counterexample being $\CN$. And the existence of functions
$\N\to\N$ which are not c.e.\ traceable shows that even first-order pointwise
finite problems may not be p.c.e.t.

The class of p.c.e.t.\ problems has some attractive algebraic properties:
\begin{itemize}
  \item If $P$ and $Q$ are p.c.e.t., then so are $P\times Q$ and $P^\ast$.
  \item If $P$ is p.c.e.t.\ and $Q\leq_\W P$, then $Q$ is p.c.e.t.
\end{itemize}
(These properties hold for problems $P$ and $Q$ on any represented spaces, as is
straightforward to show.) On the other hand,
\begin{itemize}
  \item $P\ast Q$ may not be p.c.e.t.\ even if both $P$ and $Q$ are. For
    example, $\CN$ and $\mathsf{LPO}'$ are both p.c.e.t.\ but $\CN \ast
    \mathsf{LPO}' \geq_\W \TCN \geq_\W \Ind[2]{\Q}$ \cite[Corollary~8.10]{Completion}.
  \item Neither $P^{u\ast}$ nor $P^\diamond$ need be p.c.e.t.\ if $P$ is. It
    follows from \cite[Theorem~7.2]{SV22} that $(\mathsf{LPO}')^{u \ast}
    \equiv_\W (\mathsf{LPO}')^\diamond \equiv_\W \CN'$, and
    $\CN' >_\W \TCN$ \cite[Corollary~8.14]{Completion}. (See \cite{SV22} for definitions of the
    undefined notation used here.)
\end{itemize}

Just as the fact that $\mathsf{LPO}' \vert_\W \lim$ 
(which is known in the literature, see e.g.\ \cite{Completion}) shows that
pointwise c.e.\ traceability is logically incomparable with limit computability,
the fact that $\mathsf{LPO'}$ is p.c.e.t.\ shows the concept to be strictly more general
than computability with finitely many mind changes: the latter can be
characterized by Weihrauch reducibility to $\CN$, which is p.c.e.t.\ and known
to be Weihrauch incomparable with $\mathsf{LPO}'$.

The study of p.c.e.t.\ problems has been continued by Pauly \cite{P24}, who
calls them \textit{eventually-finitely guessable} and develops the idea as part
of a family of similar notions.


\section{The countable equivalence relation}\label{sec:eq}

The countable equivalence relation $\Eqinf$ was defined in the introduction.
We identify $E=\dom \Eqinf$ with $\N\times \N$, viewing $(x,y)$ as the $y$th
element of the $x$th equivalence class. We also refer to the $x$th equivalence
class as the ``$x$th column'' of $\Eqinf$.

The proof of this theorem was obtained jointly with Linda Westrick.
\begin{thm}\thlabel{eqthm} For all $k\geq 2$, $\SRT^2_k \leq_\W \Ind[k]{\Eqinf}
  \leq_\W \RT^2_k$, but $\RT^2_2 \not\leq_\rmc \Ind[\N]{\Eqinf}$ and $\Ind[2]{\Eqinf}
  \not\leq_\W \SRT^2_\N$.
\end{thm}

\begin{figure}
\centering
\begin{tikzcd}
  \RT^2_2 \arrow[d] & \RT^2_3 \arrow[l] \arrow[d] & \RT^2_4 \arrow[l]
  \arrow[d] & \cdots \arrow[l] & \RT^2_+ \arrow[l] \arrow[d] & \RT^2_\N
  \arrow[l] \arrow[d]\\
  \Ind[2]{\Eqinf} \arrow[d] & \Ind[3]{\Eqinf} \arrow[l] \arrow[d] &
  \Ind[4]{\Eqinf} \arrow[l] \arrow[d] & \cdots \arrow[l] & \Ind{\Eqinf}
  \arrow[l] \arrow[d] & \Ind[\N]{\Eqinf} \arrow[l] \arrow[d]\\
  \SRT^2_2 & \SRT^2_3 \arrow[l] & \SRT^2_4 \arrow[l] & \cdots \arrow[l] &
  \SRT^2_+ \arrow[l] & \SRT^2_\N \arrow[l]
\end{tikzcd}
\caption{Reductions between $\SRT^2_k$, $\Ind[k]{\Eqinf}$, and
    $\RT^2_k$. An arrow from $Q$ to $P$ signifies that $P<_\W Q$. No other
    Weihrauch reductions hold other than those implied by transitivity. 
    The diagram remains true if $<_\W$ is replaced by $<_\rmc$, except that the
    reductions $\SRT^2_k \leq_\rmc \Ind[k]{\Eqinf}$ are not known to be strict,
    and that the two rightmost entries in each row collapse under $\equiv_\rmsc$.
}
    \label{eqreds}
\end{figure}
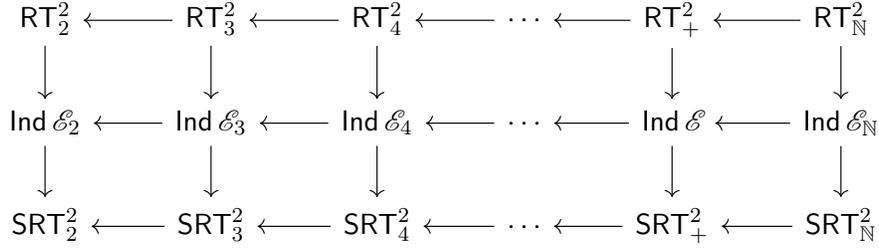

Since $\SRT^2_j \not\leq_\rmc \RT^2_k$ whenever $j>k\geq 2$ \cite{HJ16,
PateyWeakness}, we obtain Figure~\ref{eqreds}, a very similar diagram to that
shown for $\Ind\Q$ in Figure~\ref{Qreds}. To see why $\SRT^2_\N \not\leq_\W
\RT^2_+$, notice that $\SRT^2_\N$ is a fractal; the proof is along the same
lines as those of \thref{SNfractal}(ii) which implies the same for
$\Ind[\N]{\Eqinf}$. Then if there were such a reduction, we would have
$\SRT^2_\N \leq_\W \RT^2_k$ for some $k$, and in particular $\SRT^2_{k+1}
\leq_\W \RT^2_k$, a contradiction. The remaining nonreductions immediately
follow.

The rest of the section is occupied with the proof of \thref{eqthm}, which is
broken into four lemmas.
The first two establish that $\SRT^2_k \leq_\W \Ind[k]{\Eqinf} \leq_\W \RT^2_k$,
the third shows $\RT^2_2 \not\leq_\rmc \Ind[\N]{\Eqinf}$ and thus gives a
separation between $\Ind[k]{\Eqinf}$ and $\RT^2_k$, and the fourth shows
$\Ind[2]{\Eqinf} \not\leq_\W \SRT^2_\N$ and so separates $\Ind[k]{\Eqinf}$ from
$\SRT^2_k$. 

We use the notation $\{x,y\}$ for unordered pairs, and omit the extra
parentheses from $c(\{x,y\})$ and $c((x,y))$ to reduce visual clutter.

\begin{lem}\thlabel{rtdoeseq} $\Ind[k]{\Eqinf} \leq_\W \RT^2_k$.
\end{lem}
\begin{proof} Let $c\colon \N \to k$ be a coloring of $E$. Define an
  $\RT^2_k$-instance $d\colon [\N]^2\to k$ by
  \[
    d\{x,y\} = c(x,y) \quad\text{ if }x<y.
  \]
  Let $\tilde{H}$ be an infinite
  homogeneous set for $d$, and let $H$ be the $\tilde{H}$-computable set
  \[
    H = \set{ (x,y) \st x\leq y\in \tilde{H}} \subseteq \N^2.
  \]
  It is clear that $H$ is monochromatic for $c$. For each $x\in \tilde{H}$,
  $(x,y)\in H$ for all of the infinitely many $y>x$ in $\tilde{H}$, so $H$
  induces a substructure isomorphic to $\Eqinf$.
\end{proof}

\begin{lem}\thlabel{eqdoessrt} $\SRT^2_k \leq_\W \Ind[k]{\Eqinf}$.
\end{lem}
\begin{proof} Let $c\colon [\N]^2\to k$ be an instance of $\SRT^2_k$. Define
  $d\colon E \to k$ by
  \[
    d(x,y) = \begin{cases}
      c\{x,y\}, & \text{if }x<y\\
      0, & \text{otherwise.}
    \end{cases}
  \]
  Let $H$ be an $\Ind[k]{\Eqinf}$-solution for $d$. We build an $\SRT^2_k$-solution
  for $c$ as follows: start with any $(x_1,y)\in H$ with $x_1<y$.
  Then find another column $x_2>x_1$ represented in $H$ such that
  $c\{x_1,x_2\}=c\{x_1,y\}$. There is such an $x_2$
  by stability of the coloring $c$: the $d$-color of $H$ must be the same as the
  stable $c$-color of column $x_1$, so for all large enough $x$, $c\{x_1,x\} =
  c\{x_1, y\} = d(x_1,y)$. Since there are infinitely many
  equivalence classes of $E$ represented in $H$, and the equivalence class
  of $x_1$ has infinitely many points of this color, there is eventually a
  column $n$ represented in $H$ with $d(x_1,x) = d(x_1,y)$,
  and we can take $x_2$ to be this $x$. (Note $(x_1,x_2)$ need not be in
  $H$.) Next, find a column $x_3> x_2$ represented in $H$ such that
  \[
    c\{x_1,x_2\} = c\{x_1,x_3\} = c\{x_2,x_3\}.
  \]
  Again, this must be possible by stability and the existence of infinitely many
  different columns represented in $H$. Proceed in the same way, for each $i$
  finding a column $x_i$ of $H$ with $x_i>x_{i-1}$ such that $c\{x_j,x_k\} =
  c\{x_1,x_2\}$ for all $1\leq j<k\leq i$. The set $\set{x_i \st i\in\N}$ is
  then an infinite homogeneous set for $c$.
\end{proof}

It is clear that the proofs of the above lemmas do not depend on $k$ and thus
extend to show $\SRT^2_+ \leq_\W \Ind{\Eqinf} \leq_\W \RT^2_+$ and $\SRT^2_\N
\leq_\W \Ind[\N]{\Eqinf} \leq_\W \RT^2_\N$.

\begin{cor}\thlabel{eqsep1}
  For all $k\geq 2$, $\Ind[k+1]{\Eqinf} \not\leq_\rmc \Ind[k]{\Eqinf}$ and
  $\Ind{\Eqinf} \not\leq_\rmc \Ind[k]{\Eqinf}$.
\end{cor}
\begin{proof} 
  Both statements follow from the above lemmas together with the fact that
  $\SRT^2_{k+1} \not\leq_\rmc \RT^2_k$ \cite[Corollary~3.6]{PateyWeakness}.
\end{proof}

\thref{eqdoessrt} also enables us to recover the following result originally due
to Ackerman, Freer, Reimann, and Westrick \cite{Ind}:
\begin{cor} $\Eqinf$ is not computably indivisible.
\end{cor}
\begin{proof} By a result of Jockusch (see \cite[Exercise~6.35]{STT}), there is
  a computable instance of $\SRT^2_2$ with no computable solution. The same must
  be true for $\Ind[2]{\Eqinf}$ by \thref{eqdoessrt}.
\end{proof}

\begin{lem} $\RT^2_2 \not\leq_\rmc \Ind[\N]{\Eqinf}$.
\end{lem}
\begin{proof} We show that $\Ind[\N]{\Eqinf}$ admits $\Delta^0_2$ solutions,
  i.e., every computable instance has a $\Delta^0_2$ solution. As $\RT^2_2$
  does not share this property \cite[Corollary~3.2]{Jock72}, the separation
  follows.

  Let $c$ be a computable coloring of $E$ with $\ran c \subseteq k$ for some
  $k$. 
  Suppose $c$ stabilizes in infinitely many columns.
  Nonuniformly pick a color $i$ which is the stable color of infinitely many
  columns, and use $\emptyset'$ to enumerate the set $S$ of such columns by
  asking, for each $x,y$, if there is $z>y$ with $c(x,z) \neq i$, and putting
  $x\in S$ if the answer is `no'.
  As a relativization of the fact that every infinite c.e.\ set has an infinite
  computable subset, there is an infinite $\Delta^0_2$ subset of $S$, and the
  set of points of color $i$ in all such columns is a $\Delta^0_2$ solution to
  $c$.

  If instead $c$ stabilizes in only finitely many columns, nonuniformly delete
  those columns; the remaining columns each have infinitely many points of at
  least two different colors. Suppose there are colors $i_1\neq i_2$ such that
  infinitely many columns include only finitely many points of colors other than
  $i_1$ or $i_2$. The set $S$ of such columns can again be enumerated by
  $\emptyset'$ by asking, for each $x,y$, if there is $z>y$ with $c(x,z) \notin
  \{i_1,i_2\}$, and putting $x\in S$ if the answer is `no'. $S$ has an infinite
  $\Delta^0_2$ subset, and the set of points of color $i_1$ (or of $i_2$) in the
  columns in that subset is a $\Delta^0_2$ solution to $c$ because, by
  assumption, all columns of $S$ have infinitely many points both of color $i_1$
  and of color $i_2$.

  In general, let $n\geq 1$ be the least possible number of distinct colors
  $i_1, \dotsc, i_n$ such that there is an infinite set $S$ of columns each of
  which includes only finitely many points of colors other than
  $i_1,\dotsc,i_n$. That is, $n$ is such that only finitely many columns have
  infinitely many points of $n-1$ or fewer different colors. 
  Nonuniformly delete this finite set of columns, so that every column remaining
  has infinitely many points of each of at least $n$ colors, and fix a
  particular choice of $i_1, \dotsc, i_n$ as above.
  Then for this set of colors $\{i_1, \dotsc, i_n\}$, $\emptyset'$ can
  enumerate $S$ by letting $x\in S$ if we find a $y$ such that there is no $z>y$
  having $c(x,z) \notin \{i_1,\dotsc,i_n\}$. There is an infinite $\Delta^0_2$
  subset of $S$, and the set of points of color $i_j$ in the columns of this
  subset, for any $1\leq j\leq n$, is a $\Delta^0_2$ solution of $c$. \qedhere
\end{proof}

\begin{lem}\thlabel{eqsrtsep} $\Ind[2]{\Eqinf} \not\leq_\W \SRT^2_\N$.
\end{lem}
\begin{proof} $\Ind[2]{\Eqinf}$ is a total fractal and $\SRT^2_\N \equiv_\W
  \SRT^2_+ \ast \CN$ (see \thref{SNfractal} and the discussion after
  \thref{IndSN}), so if $\Ind[2]{\Eqinf} \leq_\W \SRT^2_\N$, then there is some
  $k$ with $\Ind[2]{\Eqinf} \leq_\W \SRT^2_k$. Suppose there is such a reduction
  via the functionals $\Delta$ and $\Psi$. We use lowercase Latin letters for
  instances of $\Ind[2]{\Eqinf}$ and lowercase Greek letters for their initial
  segments.

  We will find (noncomputably) a finite sequence of strings $\sigma_0 \subset
  \sigma_1 \subset \cdots$ such that for some $s$, there is an
  $\Ind[2]{\Eqinf}$-instance $c\in [\sigma_s]$ which defeats the Weihrauch
  reduction. To accomplish this, we exploit the requirement for $\Delta$ to
  produce a stable coloring in order to restrict our attention to sets of
  columns $H$ which are guaranteed to extend to an $\SRT^2_k$-solution as long
  as there is any such solution of the same color as $[H]^2$, no matter how we
  continue to extend $\sigma_s$.
  Then we use these sets as oracles to diagonalize against $\Delta$ and $\Psi$
  by adding finitely many locks on the columns of $E$ (to be explained below),
  represented by a lock function $L = \bigcup_s L_s \colon \subseteq \N \to 2$
  that will be updated finitely many times. These locks are in fact the same as
  those used in Cohen forcing with locks, but we do not formally use a notion of
  forcing here and the argument does not necessarily produce a generic $c$. 
  After we lock a column, either $\Psi$ will be forced to output a
  non-monochromatic set, or $\Delta$ will be forced to prevent $H$ from being an
  initial segment of any solution---and so will prevent there from being
  \emph{any} solution of $H$'s color. This can be carried out for every color
  used by $\Delta$, so eventually the diagonalization will succeed.

  For any string $\sigma$ and lock function $L$, let $[\sigma,L]$ be the set of
  (finite or infinite) $\tau\in [\sigma]$ such that whenever $L(x) \in \{0,1\}$
  is defined, if $(x,y)$ is newly colored by $\tau$, then 
  $\tau(x,y) = L(x)$. ``Newly colored'' simply means that $\tau$ is
  defined on $(x,y)$ while $\sigma$ is not. So $[\sigma,L]$ is the set
  of extensions of $\sigma$ which respect the lock function $L$ by only adding
  points of color $L(x)$ in column $x$. Observe that $[\sigma,L'] \subset
  [\sigma,L]$ if $L'\supset L$.

  The central claim that makes this proof work is as follows: for any $\sigma$,
  any $\SRT^2_k$-column $n$, and any lock function $L$,
  \begin{align*}
    \exists \tau\in [\sigma, L] ~ \exists i\in \N ~ \forall \rho\in
    [\tau,L] ~ ( \text{if $\Delta^\rho$ newly colors } & 
    \text{$\{n,m\}$ where $m > n$,}\\
              &\text{then $\Delta^\rho\{n,m\} = i$}).
  \end{align*}
  That is, we can pass to an extension of $\sigma$ such that among all further
  extensions respecting the lock function, $\Delta$ will never change its mind
  again about the stable color of column $n$, and in fact will never again add a
  point of any other color to it.
  We say $\tau$ $L$-forces column $n$ to be color $i$ when this happens. If the
  claim were false for some $\sigma$, $n$, and $L$, then for every $\tau\in
  [\sigma, L]$ and $i\in k$, there would be a further $\rho(\tau) \in [\tau,L]$
  such that $\Delta^{\rho(\tau)}$ adds a point not of color $i$ into column $n$.
  Letting $\tau_0 = \sigma$ and $\tau_{s+1} = \rho(\tau_s)$ for all $s$, the
  $\Ind[2]{\Eqinf}$-instance $d = \lim_s \tau_s$ is such that $\Delta^d$ has
  infinitely many points of at least two different colors in column $n$, so is
  not an $\SRT^2_k$-instance, a contradiction. A key point is that this argument
  works for any $L$ independently of $\sigma$, so that even if we change $L$
  partway through the proof we can always extend any $\sigma_s$ to a
  $\sigma_{s+1}$ which $L$-forces a fresh column to be some color. A column
  already $L$-forced remains $L'$-forced to be the same color for any $L'
  \supset L$.

  Now we describe the procedure to find $(\sigma_s)$ and $(L_s)$ in detail.
  Let $D\subseteq k$ be the set of ``already-diagonalized'' colors, starting
  with $D=\emptyset$; we will use $D$ to keep track of which colors have
  already witnessed convergence of $\Psi$ and triggered an update of the lock
  function. As long as $D \subsetneq k$,
  the colors in $D$ at stage $s$ will be exactly those of which no
  $\Delta^c$-solution can be found for any $c\in [\sigma_s,L_s]$.
  Let $L_0=\emptyset$ and let $\sigma_0$ be a string $L_0$-forcing column $0$ to
  be some color. At stage $s>0$, let $\sigma_s \in [\sigma_{s-1},L_{s-1}]$ be
  either a string which $L_{s-1}$-forces column $s$ to be some color if it had
  not yet been so forced, or an arbitrary string otherwise. 
  For each $i\notin D$, search for a pair $(x,y)$ already colored by $\sigma_s$ with
  $x\notin \dom L_{s-1}$ (i.e., with column $x$ not already locked), together
  with a finite set $H$ such that $[H]^2$ is $\Delta^{\sigma_s}$-monochromatic
  with color $i$, such that $\sigma_s$ $L_{s-1}$-forces every column in $H$ to
  be color $i$, and such that $\Psi^{\sigma_s\oplus H}(x,y) \cvg= 1$.
  Such $x$, $y$, and $H$ must eventually be found for some $i$ and $s$ if
  $\Psi$ is to compute an $\Ind[2]{\Eqinf}$-solution, and we can take $x\notin
  \dom L_{s-1}$ because $\Psi$ will output elements of infinitely many columns
  of $E$ whereas $\dom L_{s-1}$ is finite. 
  That we can also take $i\notin D$ is justified below.
  If $x$ and $y$ are found at stage $s$, diagonalize by setting $L_s(x) = 1-
  \sigma_s(x,y)$; that is, we lock column $x$ to be the opposite color that
  $\Psi$ already committed to there. End stage $s$ by putting $i\in D$. If
  instead no such $x$, $y$, or $H$ are found for any $i$, let $L_s = L_{s-1}$
  and end the stage.

  If $H$ is as in the previous paragraph, then $H$ is a valid initial segment of
  some solution to $\Delta^d$ for \emph{any} $d\in [\sigma_s,L_{s-1}]$ for which
  infinitely many $\SRT^2_k$-columns stabilize to $i$, and in particular for any
  such $d\in [\sigma_s,L_s]$. This is because starting from any
  $\Delta^d$-solution $G$, one can truncate $G$ to $\hat{G}$ by selecting only
  the columns whose indices are higher than the point at which any of the
  columns in $H$ stabilize, producing an infinite $\hat{G}$ with $[H \cup
  \hat{G}]^2$ monochromatic. Hence either $H$ extends to a solution of color $i$
  (for some $d$ as above), or there is no solution of color $i$ (for any $d$ as
  above). In the first case, $d$ would witness the failure of the Weihrauch
  reduction, because if $H$ extends to a $\Delta^d$-solution $K$, then
  $\Psi^{d\oplus K}$ outputs a point in
  $E$ whose color is shared by only finitely many other points in the same
  column. If such a $d$ exists then the procedure ends at this stage. In
  particular, this must happen if $D=k$, because then some $H$, $x$, and $y$
  were already found for every $i\in k$, and for any $d\in [\sigma_s, L_s]$, at
  least one of these $H$s must extend to a valid solution of $\Delta^d$.

  Otherwise, for this $i$ and indeed for all $i\in D$, there is no $d\in
  [\sigma_s, L_s]$ with $\Delta^d$ having a solution of color $i$. This means
  that only finitely many columns of any such $\Delta^d$ stabilize to a color in
  $D$, and it must still be possible to $L_s$-force infinitely many columns to
  be a color not in $D$. Hence we can continue the procedure to extend
  $\sigma_s$ and search for a suitable $H$ of some color $i\notin D$.
  This completes the proof, since the procedure will eventually diagonalize
  against every color of which it is still possible for there to be an
  $\SRT^2_k$-solution. \qedhere
\end{proof}


\section{Further directions}\label{sec:moredirs}

Solutions of $\Ind{\Q}$ are rather nebulous: for example, if $x$ and $y$ are two
elements of a solution $H$, one can delete the whole interval $[x,y]$ from $H$
and still obtain a solution. The seeming weakness of $\Ind{\Q}$ may be a
consequence of this nebulosity, and it is unclear how much power 
it derives from the fact that it outputs a second-order object, i.e., a
set of rationals. 
An investigation of these properties can be put on precise footing with the
following notions introduced in \cite{DSY} and \cite{GPV21}, respectively:
\begin{itemize}
  \item The \textit{first-order part} of a problem $P$, $^1P$, is the strongest
    first-order problem Weihrauch reducible to $P$.
  \item The \textit{deterministic part} of $P$, $\Det P$, is the strongest
    problem Weihrauch reducible to $P$ for which every instance has a unique
    solution.
\end{itemize}
\begin{qu} How strong exactly are $^1\Ind{\Q}$ and $\Det(\Ind{\Q})$?
\end{qu}

The first-order part of the related tree pigeonhole principle $\mathsf{TT}^1_+$,
which was discussed in Section~\ref{sec:prior},
has been studied in recent work by Dzhafarov, Solomon, and Valenti \cite{DSV}.
They showed that none of the problems $\mathsf{TT}^1_k$, $\mathsf{TT}^1_+$, or
$\mathsf{TT}^1_\N$ are Weihrauch equivalent to any first-order problem, and
indeed $^1 \mathsf{TT}^1_\N \equiv_\W \RT^1_\N$. But $^1\mathsf{TT}^1_k$
turns out to be strictly stronger than $\RT^1_k$ for all $k\geq 2$; whether $^1
\mathsf{TT}^1_k$ and $^1 \mathsf{TT}^1_+$ have precise characterizations in
terms of known Weihrauch degrees is open, and the same applies to the
first-order parts of versions of $\Ind{\Q}$.
Regarding $\Det (\Ind{\Q})$, Manlio
Valenti has pointed out to the author that since $\Det (\RT^1_k) \equiv_\W \lim_k$ and
$\Det (\TCN^\ast) \equiv_\W \CN \equiv_\W \lim_\N$, where $\lim_k$ maps an
eventually constant element of $k^\N$ to its limit and similarly for $\lim_\N$,
we must have $\lim_k <_\W \Det (\Ind{\Q}) <_\W \lim_\N$ for all $k$. (The second
reduction is strict by \thref{QperpCN}.) Precise characterizations of the
degrees of $\Det(\mathsf{TT}^1_k)$, $\Det(\mathsf{TT}^1_+)$,
$\Det(\Ind[k]{\Q})$, and $\Det(\Ind{\Q})$ have not been established. 
On the other hand, $\Det(\Ind[\N]{\Q}) \equiv_\W \Det( \mathsf{TT}^1_\N)
\equiv_\W \CN$ since $\CN$ is below both $\Ind[\N]{\Q}$ and $\mathsf{TT}^1_\N$
and, using \cite[Theorem~3.9]{GPV21}, $\Det(\TCN^\ast \ast \CN) \leq_\W
\Det(\TCN^\ast) \ast \CN \equiv_\W \CN \ast \CN \equiv_\W \CN$.

One consequence of the above facts is that the
deterministic part of $\Ind{\Q}$ is p.c.e.t., and so one cannot replace
$\Ind[2]{\Q}$ with $\Det (\Ind[2]{\Q})$ in \thref{pcet}. 
Pending clarification of the exact relationship between $\Ind{\Q}$ and
$\mathsf{TT}^1_+$, the same should apply to $^1\Ind[2]{\Q}$.

Turning to $\Eqinf$, \thref{eqthm} is intriguing given the great interest in
$\RT^2_2$ in reverse mathematics since the 1970s. 
Cholak, Jockusch, and Slaman in \cite{CJS} proved that $\RT^2_2$ can be
logically decomposed into the conjunction of $\SRT^2_2$ and $\COH$,
where $\COH$ is the so-called cohesive principle which states that for any
sequence $(R_i)$ of subsets of $\N$, there is a set which for each $i$ is up to
finite error a subset of either $R_i$ or its complement.
This work sparked a major investigation into the logical separation of
$\SRT^2_2$ from $\COH$ which was completed a few years ago by Monin and Patey
\cite{MoninPatey}, who showed that $\SRT^2_2$ does not imply $\COH$ even in
$\omega$-models of $\mathsf{RCA}_0$, the converse nonimplication having been
established earlier in \cite{HJKLS}.
One is led to wonder what relationship $\Ind[2]{\Eqinf}$ has to
$\COH$: clearly $\COH$ cannot imply $\Ind[2]{\Eqinf}$ in $\omega$-models, but
whether $\COH$ is reducible in any sense to $\Ind[2]{\Eqinf}$ is open, as is the
question of the separation between $\Ind[k]{\Eqinf}$ and $\SRT^2_k$ under
$\leq_\rmc$.
We conjecture that
\begin{conj}\thlabel{eqconj} $\COH \not\leq_\rmc \Ind[\N]{\Eqinf}$ and
  $\Ind[2]{\Eqinf} \not\leq_\rmc \SRT^2_\N$, or at least $\Ind[k]{\Eqinf}
  \not\leq_\rmc \SRT^2_k$.
\end{conj}

Of course, if either statement is false, then the truth of the other would
immediately follow.
We briefly mention a way to view $\Ind[k]{\Eqinf}$ in terms of $\RT^1_k$ which
might suggest a route to the resolution of both conjectures.
Several authors have observed that $\COH \equiv_\W \widehat{(\RT^1_2)^\FE}$,
where $P^\FE$ is the ``finite error'' version of $P$ with $\dom P^\FE = \dom P$
and $x\in P^\FE(p)$ iff there is a $y\in P(p)$ with $x$ and $y$ differing only
on a finite set. (See for instance \cite[Corollary~8.4.15]{DM}.)
The next definition was inspired by the principle $\mathsf{RCOH}$, a
weakening of $\COH$ introduced in \cite{CDHP} which corresponds to ``solving
infinitely many columns'' of an $\widehat{(\RT^1_2)^\FE}$-instance.

\begin{defi} The \bdef{weak parallelization} of a problem $P$ is the problem
  $\widetilde{P}$ such that $\dom \widetilde{P} = \dom \widehat{P}$, and where
  the solutions to the instance $\ang{p_0,p_1,\dotsc}$ are all sets of the form
  \[
    \bigcup_{n\in A} \set{ \ang{n,x} \st x\in P(p_n)}, 
  \]
  where $A$ is an infinite subset of $\N$.
\end{defi}

In other words, $\widetilde{P}$ picks infinitely many parallel instances of $P$
to solve out of a given instance of $\widehat{P}$. Thus $\Ind[k]{\Eqinf}$ is an a
priori stronger variant of $\widetilde{\RT^1_k}$ in which the solutions of the
parallel $\RT^1_k$-instances represented in an $\widetilde{\RT^1_k}$-solution are
all of the same color. However, by nonuniformly picking a color shared by
infinitely many columns of a solution of $\widetilde{\RT^1_k}$, it is not hard
to see that $\Ind[k]{\Eqinf} \equiv_\rmc \widetilde{\RT^1_k}$. Moreover, $\SRT^2_k
\equiv_\rmc \widetilde{\SRT^1_k}$ by additionally computing a homogeneous set
using a standard argument as in the proof of \thref{eqdoessrt}. Then
\thref{eqconj} can be rephrased in these terms as
\begin{conj}
  $\widehat{(\RT^1_2)^\FE} \not\leq_\rmc \widetilde{\RT^1_2}$ and
  $\widetilde{\RT^1_k} \not\leq_\rmc \widetilde{\SRT^1_k}$.
\end{conj}

Finally, there are of course many indivisible structures other than $\Q$ and
$\Eqinf$ which could be investigated along similar lines as in the present work,
such as the Rado graph or nonscattered linear orders.
As a followup to \thref{nQm}, it would be interesting to obtain a
characterization of the linear orders $\mathcal{L}$ such that $\Ind{\mathcal{L}}
\equiv_\W \Ind{\Q}$.
It might also be fruitful to continue the study of general
properties of indivisibility problems initiated in Section~\ref{sec:rmk}.
Arguably one of the most basic questions we leave open is the following:
\begin{qu}
  Is it the case that $\Ind[k+1]{\calS} \not\leq_\W \Ind[k]{\calS}$ for every
  indivisible structure $\calS$?
\end{qu}


\bibliographystyle{alphaurl}
\bibliography{indsources}

\end{document}